\documentclass[12pt]{article}
\usepackage[final]{epsfig}
\usepackage{graphics}
\usepackage{amsmath}
\usepackage{amsfonts}
\usepackage{latexsym}
\usepackage{amssymb}
\usepackage{graphicx}
\usepackage{url}
\usepackage{epstopdf}
\usepackage{hyperref}
\usepackage{color}
\usepackage{marginnote}
\usepackage{a4wide}

\newtheorem{lemma}{Lemma}[section]
\newtheorem{proposition}[lemma]{Proposition}
\newtheorem{remark}[lemma]{Remark}

\newtheorem{theorem}{Theorem}

\newtheorem{conjecture}[lemma]{Conjecture}

\newcommand{\eps}{{\varepsilon}}
\newcommand{\proofend}{$\Box$\bigskip}

\newcommand{\R}{{\mathbb R}}

\def\proof{\paragraph{Proof.}}

\newcommand{\T}{\stackrel{t}{\sim}}

\begin{document}

\title{A family of maps and a vector field on plane polygons}

\author{Maxim Arnold\footnote{
Department of Mathematics,
University of Texas Dallas,
Richardson, TX 75080, USA;
Maxim.Arnold@utdallas.edu}
 \and 
 Lael Costa\footnote{
Department of Mathematics,
Pennsylvania State University,
University Park, PA 16802,
USA;
lael.costa@psu.edu}
\and
 Serge Tabachnikov\footnote{
Department of Mathematics,
Pennsylvania State University,
University Park, PA 16802,
USA;
tabachni@math.psu.edu}
}  

\date{\today}

\maketitle

\begin{abstract}
We study, theoretically and experimentally, a 1-parameter family of transformations and their limiting vector field on the space of plane polygons. These transformations are discrete analogs of completely integrable transformation on closed plane curves, known as the bicycle correspondence, that is a geometric realization of the B\"acklund transformation of the planar filament equation. For odd-gons, we construct a symplectic form 
on the quotient space by parallel translations and show that the transformations are symplectic, and  the vector field is Hamiltonian. 
In the case of triangles, we prove complete integrability of the respective vector field and provide evidence for the conjecture that the transformations are  integrable as well.
\end{abstract}

\section{Introduction} \label{sect:intro}
The motivation for this work comes from the recent studies of a simple model of bicycle kinematics. In this model, bicycle is represented by an oriented segment of fixed length that can move in the plane in such a way that the velocity of its rear end is always aligned with the segment (the rear wheel of the bicycle is fixed on the frame). In this model, the rear track determines the front track  as the locus of the end points of the positive tangent segments of a fixed length. See, e.g., \cite{FLT}.

A closed rear track defines, depending on its orientation, two front tracks, and one says that they are in the bicycle correspondence. This relation on closed plane curves is completely integrable: it is a geometric realization of the B\"acklund transformation of the (planar)  filament equation. See \cite{BLPT,Ta}.

The bicycle correspondence was discretized in \cite{TT} where it is defined as a relation on plane polygons. In the present paper we present and study a different discretization of the relation between rear and front bicycle tracks and of the bicycle correspondence. 

\begin{figure}[ht]
\centering
\includegraphics[width=.45\textwidth]{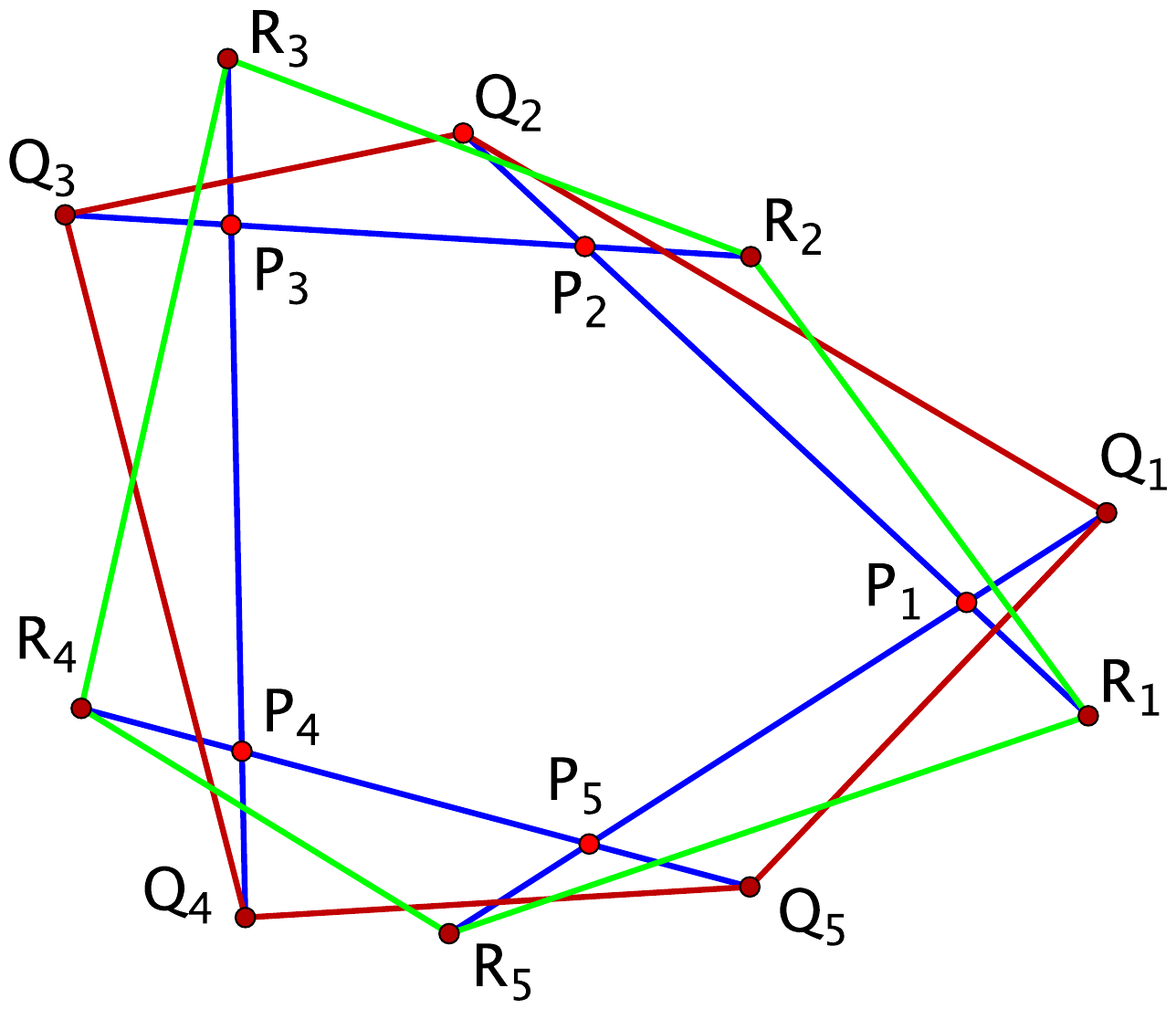}\qquad
\includegraphics[width=.45\textwidth]{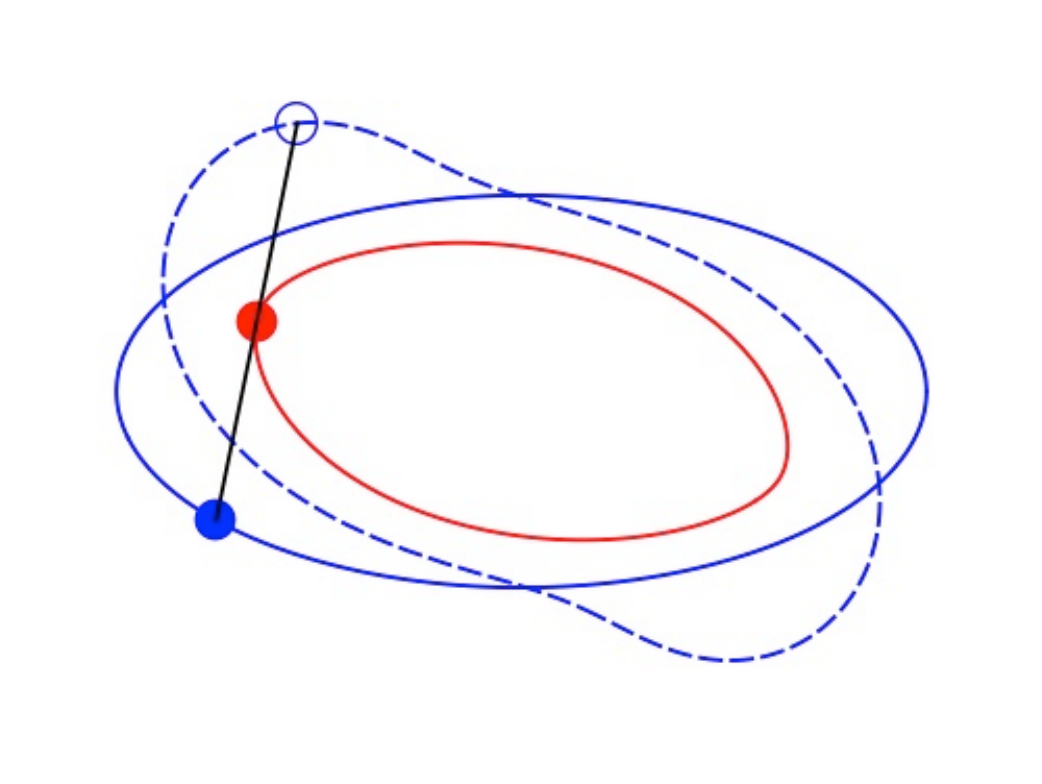}
\caption{Left: ${\bf Q} = q_t({\bf P}), {\bf R} = r_t({\bf P})$. Right: two curves in the bicycle correspondence.}	
\label{map}
\end{figure}

Fix $t>0$ (``the bicycle frame"), and let ${\bf P}=(P_1 P_2 \ldots P_n)$ be a plane $n$-gon (``the  rear track"). We define two new $n$-gons, ${\bf Q}$ and ${\bf R}$ (``the front tracks"):
$$
Q_i = P_i + t\frac{P_i - P_{i-1}}{|P_i - P_{i-1}|}, \ R_i = P_i + t\frac{P_i - P_{i+1}}{|P_i - P_{i+1}|}, \ i=1,\ldots, n,
$$
see Figure \ref{map}. Here and elsewhere the indices are understood cyclically. We have defined two 1-parameter families of maps
$ q_t ({\bf P}) = {\bf Q},  r_t({\bf P}) = {\bf R}.$

This defines a relation between the polygons ${\bf Q}$ and ${\bf R}$ that we denote by ${\bf Q} \T {\bf R}$. This relation 
 is a discrete analog of the above described bicycle correspondence. 

Let us briefly describe the contents of the paper.

In Section \ref{sect:map} we  show that, for sufficiently small values of $t$, the relation $\T$ is actually a map, $f_t$. In the limit $t\to 0$, one obtains a vector field, $\xi$.

Section \ref{sect:odd} concerns $n$-gons with $n$ odd. We consider the quotient space of such polygons by parallel translations and define a symplectic structure on this moduli space. The relation $\T$ descends on this moduli space as a symplectic relation, and the vector field $\xi$ is Hamiltonian with the perimeter of a polygon as the Hamiltonian function. The relations $\T$ have an integral, a certain quadratic function of the coordinates of the vertices of a polygon that we call the algebraic multi-area. The respective Hamiltonian vector field is an infinitesimal rotation. 

In Section \ref{sect:tri} we study the case of triangles in a greater detail. The vector field $\xi$ is Liouville integrable: its Poisson commuting integrals are the perimeter and the area of a triangle. Using the area function, one performs symplectic reduction, arriving at the 2-dimensional space of shapes of triangles (it is itself a triangle). We calculate the resulting area form on this shape space and the integral of the vector field $\xi$ that descends on this space. This field has a single zero, corresponding to equilateral triangles, and all its trajectories are closed. We also show that two triangles in the relation $\T$ are Poncelet: they are inscribed in a conic and circumscribed about a conic.

Section \ref{sect:even} concerns even-gons. Associated with such a polygon is the alternating sum of its vertices, a well defined vector. We show that this vector is invariant under $\T$, thus providing two integrals of this relation. The perimeter of a polygon is invariant under the vector field $\xi$, just as in the case of odd-gons.

The last Section \ref{sect:exp} contains open questions, conjectures, and
various experimental results. Our experimental tool is a piece of software developed by
the second-named author and available at
\url{https://lael.dev/#/evasion}. The code may be found at
\url{https://www.github.com/Lael/MathDemos}.

\bigskip

{\bf Acknowledgements}. The authors participated in the special program ``Mathematical Billiards: at the Crossroads of Dynamics, Geometry, Analysis, and Mathematical Physics" at Simons Center for Geometry and Physics where this work was started. We are grateful to the center for its hospitality and stimulating atmosphere. We thank A. Akopyan and V. Zharnitsky for useful discussions. ST and LC were supported by NSF grant DMS-2005444.

\section{A family of maps and a vector field} \label{sect:map}

We have defined a family of relations (that is, a multi-valued maps) ${\bf Q} \T {\bf R}$. 
We would like to think about this relation as a map, at least for sufficiently small values of $t$. For that, one needs to be able to construct the inverse maps of $q_t$ and $r_t$. Once it is done, one defines $f_t = q_t \circ r_t^{-1}$, so  that $f_t ({\bf R}) = {\bf Q}$. Likewise, $f_t^{-1} = r_t \circ q_t^{-1}$.

Let us describe the construction of $q_t^{-1}$ (for $r_t^{-1}$, the construction is similar). Given an $n$-gon ${\bf Q}$, assume that $t$ is less than half of the length of every side of ${\bf Q}$.

\begin{lemma} \label{lm:constr}
There exists a unique $n$-gon ${\bf P}$ such that $q_t ({\bf P}) = {\bf Q}$.
\end{lemma}

\begin{figure}[ht]
\centering
\includegraphics[width=.45\textwidth]{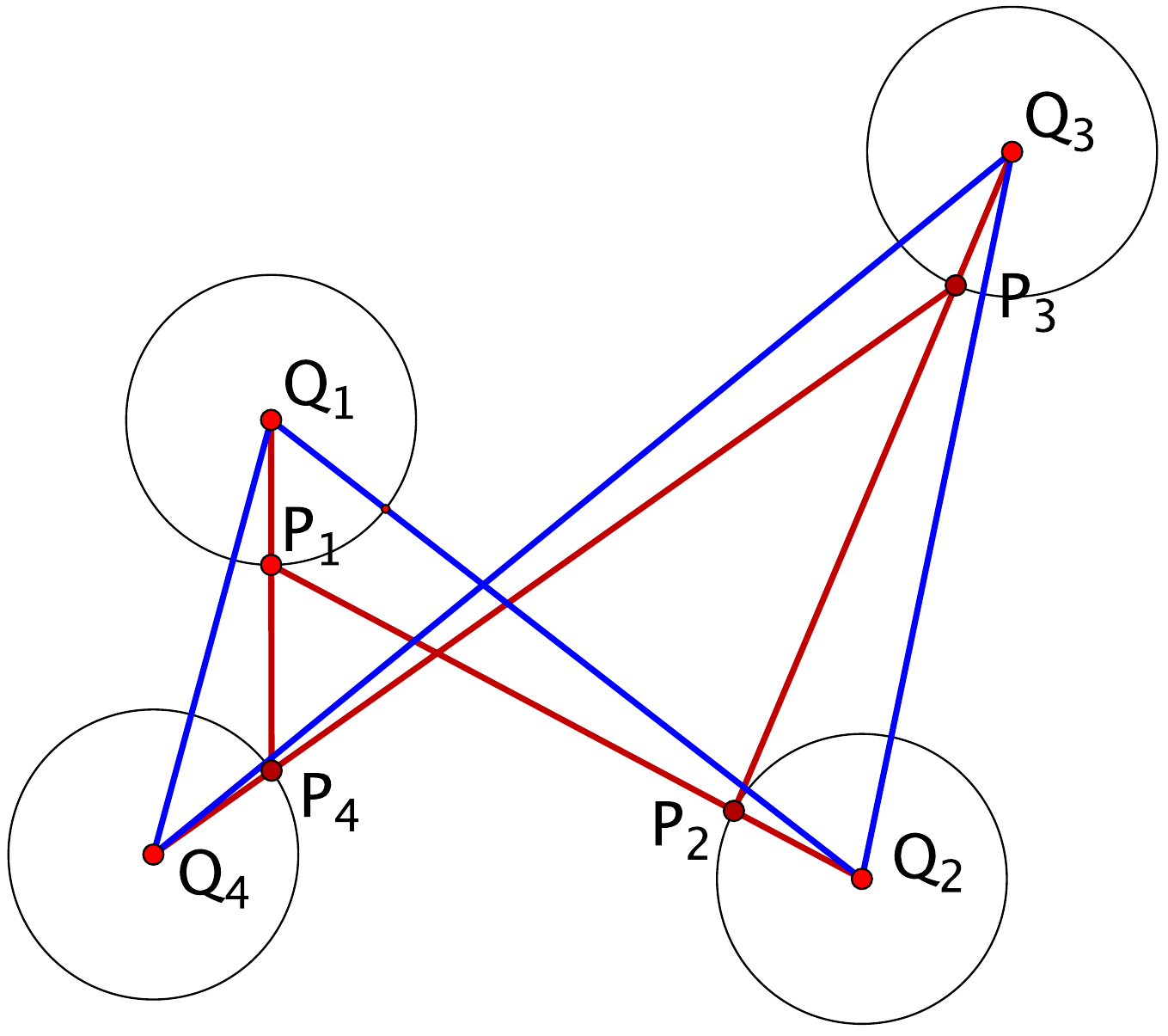}
\caption{Constructing $q_t^{-1}$.}	
\label{inverse}
\end{figure}

\proof
Consider circles of radii $t$ centered at every vertex of ${\bf Q}$. See Figure \ref{inverse}. 

Choose a test point $P_1$ on the circle centered at $Q_1$. Since $2t< |Q_1 Q_2|$, one has 
$|P_1 Q_2| > |Q_1 Q_2| - t > t$. Hence the segment $P_1 Q_2$  intersects the circle centered at $Q_2$; let $P_2$ be this intersection point. Next, let  $P_3$ be the intersection point of the segment $P_2 Q_3$ with the circle centered at $Q_3$, and so on. The $(n+1)$st point $P_{n+1}$ again lies on the circle centered at $Q_1$. 

Note that the central projection of the exterior of a circle on the circle is distance decreasing. Therefore each step of the above described construction is a contracting map from each circle to the next one.  It follows that  there is a unique choice of $P_1$ such that $P_{n+1}=P_1$. This gives an $n$-gon ${\bf P}$ such that $q_t ({\bf P})={\bf Q}$. 
\proofend

\begin{figure}[ht]
\centering
\includegraphics[width=.4\textwidth]{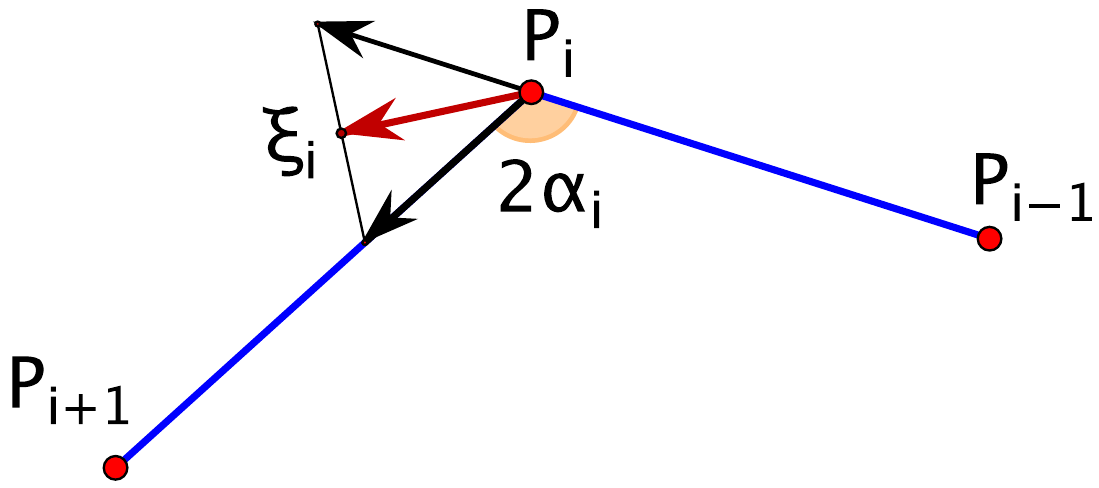}
\caption{The vector field $\xi$.}	
\label{field}
\end{figure}

We also define a vector field $\xi$ on the space of $n$-gons that is a limiting case of the maps $f_t$ as $t\to 0$. This vector field $\xi=(\xi_1,\xi_2,\ldots,\xi_n)$ is given by
$$
\xi_i = \frac{1}{2} \left(\frac{P_i - P_{i-1}}{|P_i - P_{i-1}|} - \frac{P_i - P_{i+1}}{|P_i - P_{i+1}|}\right),
$$
see Figure \ref{field}.

The vector $\xi_i$ has the direction of the bisector of the exterior angle of the polygon at vertex $P_i$, and $|\xi_i|=\sin \alpha_i$, where $2\alpha_i$ is the interior angle of the polygon at vertex $P_i$.

\section{Odd-gons} \label{sect:odd}

Let $n$ be odd, and ${\bf P}$ be an $n$-gon. In Cartesian coordinates, one has  $P_i=(x_i,y_i)$. Thus $(x_1,\ldots,x_n,y_1,\ldots,y_n)$ are coordinates on the space of $n$-gons, identified with $\R^{2n}$.

Define a closed differential 2-form on the space of $n$-gons:
$$
\omega=\sum_{1\le i < j \le n} (-1)^{i+j-1} (dx_i \wedge dx_j + dy_i\wedge dy_j).
$$

\begin{lemma} \label{lm:form}
The 2-form $\omega$ has a 2-dimensional kernel, generated by parallel translations acting diagonally on polygons.
\end{lemma}

\proof
The matrix of $\omega$ is a direct sum of two identical circulant $n\times n$ matrices (in $x$ and in $y$ variables) whose first
row is $(0,1,-1,1,-1,\ldots,-1)$. The eigenvalues of such circulant matrix are $h(\lambda_k),\ k=0,1,\ldots,n-1$, where 
$$
h(t)=x-x^2+x^3-\ldots + x^{n-2}-x^{n-1}
$$
 and $\lambda_k=e^{2\pi ik/n}$, a root of unity, see \cite{Da}. One has
$$
h(x)=x\frac{1-x^{n-1}}{1+x},
$$
and $h(\lambda_k)=0$ for $k=0$ only. 

The respective eigenvector is $(1,\ldots,1)$; it corresponds to parallel translations in the $x$-direction. Likewise for the $y$-component.  
\proofend

Thus $\omega$ becomes non-degenerate on space of $n$-gons factorized by parallel translations. We will work with this moduli space; denote it by ${\mathcal T}$ or, when one needs to specify $n$, by ${\mathcal T}_n$. Abusing notation, we use the same symbols, $\xi$ and $f_t$, for the induced vector field and maps on ${\mathcal T}$. 

In terms of the $(x,y)$-coordinates, our vector field is given by the formula
$$
\xi = \sum \left(\frac{x_i - x_{i-1}}{a_i} + \frac{x_{i+1}-x_i}{a_{i+1}}\right) \partial x_i +
\left(\frac{y_i - y_{i-1}}{a_i} + \frac{y_{i+1}-y_i}{a_{i+1}}\right) \partial y_i,
$$
where $a_i=|P_i-P_{i-1}|$ is the side length, and the sums are cyclic.

Let ${\mathcal P}=\sum a_i$ be the perimeter function on ${\mathcal T}$.

\begin{lemma} \label{lm: Ham}
The field $\xi$ is Hamiltonian with the Hamiltonian function ${\mathcal P}$.
\end{lemma}

\proof
One has
$$
{\mathcal P} = \sum \sqrt{(x_{i+1}-x_i)^2+(y_{i+1}-y_i)^2},
$$
and 
$$
d {\mathcal P} = \sum \left(\frac{x_i - x_{i-1}}{a_i} - \frac{x_{i+1}-x_i}{a_{i+1}}\right) dx_i +
\left(\frac{y_i - y_{i-1}}{a_i} - \frac{y_{i+1}-y_i}{a_{i+1}}\right) dy_i .
$$
One checks that this equals $i_\xi \omega$.
\proofend

Define another function on ${\mathcal T}$, an algebraic multi-area:
$$
{\mathcal A}({\bf P})= \sum_{1\le i < j \le n} (-1)^{i+j} \det(P_i,P_j) = \sum_{k=1}^{(n-1)/2)} (-1)^k \sum_{i=1}^n \det(P_i,P_{i+k}).
$$

Let $\eta = \sum x_i \partial y_i - y_i \partial x_i$ be the rotation vector field acting on polygons diagonally. 

\begin{lemma}  \label{lm:rot}
One has $d {\mathcal A} = i_\eta \omega$.
\end{lemma}

\proof
It is a direct calculation, and we omit it.
\proofend

Since perimeter is invariant under rotations, the functions ${\mathcal A}$ and ${\mathcal P}$ Poisson commute, and ${\mathcal A}$ is an integral of the vector field $\xi$. 

Let us show that the algebraic multi-area is an integral of the relation $\T$.

\begin{lemma} \label{lm:Int}
If ${\bf Q} \T {\bf R}$ then ${\mathcal A}({\bf Q})= {\mathcal A}({\bf R})$.
\end{lemma}

\proof
Let $e_i$ be the unit vector from $P_{i-1}$ to $P_i$. Then
$
Q_i=P_i-te_{i+1},\ R_i=P_i+te_i.
$
One has
$$
\begin{aligned}
&{\mathcal A}({\bf Q})=\sum_{k=1}^m (-1)^k \sum_{i=1}^n \det(Q_i,Q_{i+k}) = \sum_{k=1}^m (-1)^k \sum_{i=1}^n \det(P_i-te_{i+1},P_{i+k}-te_{i+k+1})=\\
&{\mathcal A}({\bf P}) + t \sum_{k=1}^m (-1)^k \sum_{i=1}^n [-\det(P_i,e_{i+k+1})+\det(P_{i+k},e_{i+1})] +t^2 \sum_{k=1}^m (-1)^k \sum_{i=1}^n \det(e_{i+1},e_{i+k+1}).
\end{aligned}
$$
Likewise,
$$
{\mathcal A}({\bf R})={\mathcal A}({\bf P}) + 
t \sum_{k=1}^m (-1)^k \sum_{i=1}^n [-\det(P_i,e_{i+k})+\det(P_{i+k},e_{i})] +
t^2 \sum_{k=1}^m (-1)^k \sum_{i=1}^n \det(e_{i},e_{i+k}).
$$
Clearly, the terms of degree zero and two in $t$ are the same in ${\mathcal A}({\bf Q})$ and ${\mathcal A}({\bf R})$. It remains to show that the linear in $t$ terms are  equal as well.

Since the sums in $i$ are cyclic, it suffices to consider the terms of the form $\det(P_1,e_j)$ in both sums. Denote the the linear in $t$ part of
these terms in  ${\mathcal A}({\bf Q})$ by $q_j$, and in ${\mathcal A}({\bf R})$ by $r_j$. Then 
$$
q_1=-1,q_2=0,q_3=1,q_4=-1,q_5=1,\ldots,q_{m+1}=(-1)^m,\ldots, q_n=1,
$$
and
$$
r_1=0,r_2=-1,r_3=1, r_4=-1, r_r=1,\ldots,q_{m+1}=(-1)^m,\ldots, q_n=1.
$$
Therefore  
$$
\begin{aligned}
&{\mathcal A}({\bf Q }) - {\mathcal A}({\bf R})=\sum_i  [\det(P_i,e_{i+1})-\det(P_{i},e_i)] =\\
&\sum_i  [\det(P_{i-1},e_{i})-\det(P_{i},e_i)] = \sum_i \det(P_{i-1}-P_i,e_i) =0,
\end{aligned}
$$
the last equality holds because $e_i$ is collinear with $P_i-P_{i-1}$.
\proofend

Furthermore, we have

\begin{lemma} \label{lm:Inv}
The relation $\T$ is symplectic.
\end{lemma}

\proof
One needs to check that $q_t^* (\omega)=r_t^* (\omega)$. As before, both sides are quadratic polynomials in $t$.

Let $e_i=(u_i,v_i)$. One needs to show that quadratic polynomials in $t$ are equal, and it is clear that the free and quadratic terms are. As in the proof of Lemma \ref{lm:Int}, one needs to focus on the linear terms. That is, one needs to compare
$$
\sum_{k=1}^m (-1)^k \sum_{i=1}^n [ - dx_i \wedge du_{i+k+1} + dx_{i+k} \wedge du_{i+1}
- dy_i \wedge dv_{i+k+1} + dy_{i+k} \wedge dv_{i+1}]
$$
with
$$
\sum_{k=1}^m (-1)^k \sum_{i=1}^n [dx_i \wedge du_{i+k} - dx_{i+k} \wedge du_i +
dy_i \wedge dv_{i+k} - dy_{i+k} \wedge dv_i].
$$

We need to show that these two sums are equal. The same analysis as in the proof of Lemma \ref{lm:Int} reduces this to showing that
\begin{equation} \label{eq:sumw}
\sum_i  [(dx_{i+1}-dx_i) \wedge du_{i+1} + (dy_{i+1}-dy_i) \wedge dv_{i+1}] =0.  
\end{equation}

Recall that $a_i$ denote the side lengths. One has
$$
\begin{aligned}
&du_{i+1} = \frac{(y_{i+1}-y_i)}{a_{i+1}^3} [(y_{i+1}-y_i) (dx_{i+1}-dx_i) - (x_{i+1}-x_i) (dy_{i+1}-dy_i)],\\
&dv_{i+1} = -\frac{(x_{i+1}-x_i)}{a_{i+1}^3} [(y_{i+1}-y_i) (dx_{i+1}-dx_i) - (x_{i+1}-x_i) (dy_{i+1}-dy_i)].
\end{aligned}
$$
Substitute to (\ref{eq:sumw}): all the terms cancel, and the sum equals zero, as needed.
\proofend

Lemmas \ref{lm:rot}, \ref{lm:Int}, and \ref{lm:Inv} make it possible to perform symplectic reduction. Fix a non-zero (hence regular) value of the function ${\mathcal A}$ and factorize by the rotations. The resulting quotient space can be identified with the quotient space of the space of $n$-gons by orientation preserving homotheties; we call this symplectic $2n-4$-dimensional space the {\it shape space} and denote it by ${\mathcal S}$. The relation $\T$ descends to ${\mathcal S}$ as a symplectic relation and the vector field $\xi$ as a Hamiltonian vector field.
We continue to abuse notation and use the same symbols for the resulting vector field and the maps on ${\mathcal S}$.

In conclusion of this section, let us explain the origin of the double sums occuring in the formulas for $\omega$ and ${\mathcal A}$. 

Let ${\bf P}$ be an $n$-gon, and recall that $n$ is odd. 

\begin{lemma} \label{lm:out}
There exists a unique $n$-gon ${\bf S}$ such that the vertices of ${\bf P}$ are the midpoints of the sides of ${\bf S}$.
Then $2 {\mathcal A} ({\bf P})$ equals the (signed) area of ${\bf S}$. 
\end{lemma}

\proof
One needs to solve the system of linear equations
$$
P_i = \frac{S_i+S_{i+1}}{2},\ i=1,\ldots,n.
$$
Since $n$ is odd, the system has a unique solution
$$
S_i = \sum_{j=0}^{n-1} (-1)^j P_{i+j},\ i=1,\ldots,n,
$$
as claimed.

Using this formula and collecting terms in the formula for the area
$$
\frac{1}{2} \sum \det (S_i,S_{i+1})
$$
yields $2 {\mathcal A} ({\bf P})$.
\proofend

Let $S_i=(a_i,b_i)$, and let $\Omega = \sum (da_i \wedge da_{i+1} + db_i \wedge db_{i+1})$. The same calculation as in the previous lemma yield

\begin{lemma} \label{lm:Symp}
One has $\Omega = 4 \omega$.
\end{lemma}

Finally, we interpret the vector field $\xi$ in terms of the polygons ${\bf S}$.  

Let $\eta$ be the vector field where vertex $S_i$ moves with the unit speed in the direction of the diagonal $S_{i-1} S_{i+1}$:
$$
\eta_i = \frac{S_{i+1} - S_{i-1}}{|S_{i+1} - S_{i-1}|}.
$$

\begin{lemma} \label{lm:fields}
One has $\eta=2 \xi$.
\end{lemma}

\proof
Since 
$$
P_i = \frac{S_i+S_{i+1}}{2},\ S_{i+1}-S_{i-1}=2(P_i-P_{i-1}),
$$
 the velocity of $P_i$ is 
$$
\frac{1}{2} \left(\frac{S_{i+1} - S_{i-1}}{|S_{i+1} - S_{i-1}|} + \frac{S_{i+2} - S_{i}}{|S_{i+2} - S_{i}|}\right)=
\frac{P_{i} - P_{i-1}}{|P_{i} - P_{i-1}|} + \frac{P_{i+1} - P_{i}}{|P_{i+1} - P_{i}|} = 2 \xi_i,
$$
as claimed.
\proofend

Note that the perimeter of $P$ equals one half of the perimeter of the ``star-shaped" polygon $S_1 S_3 S_5\ldots$, therefore the perimeter of this star-shaped polygon is preserved by the vector field $\eta$.

\section{Triangles} \label{sect:tri}

In this section, we specialize the results of the previous section to triangles and investigate this case in more detail. 

The  space of triangles modulo parallel translations is a 4-dimensional symplectic space ${\mathcal T}$. On this space we have the Hamiltonian vector field $\xi$ and a family of the symplectic relations $\T$.
The results of Section \ref{sect:odd} imply the following theorem.

\begin{theorem} \label{thm:tri}
The vector field $\xi$ is completely integrable: it has Poisson commuting integrals, the perimeter and the area of a triangle. 
The relation $\T$ has a conserved quantity, independent of $t$, the area of a triangle.
\end{theorem}

As described in the preceding section, we perform symplectic reduction. The resulting shape space ${\mathcal S}$ consists of isometry equivalence classes of triangles with a fixed (say, unit) area or, equivalently, of the orientation preserving homothety equivalence classes of triangles.

The shape of a triangle is determined by its angles $\alpha, \beta, \gamma$ satisfying $\alpha+\beta+\gamma=\pi$. Therefore the shape space is itself a triangle:
$$
{\mathcal S}= \{(\alpha,\beta,\gamma)\ | \ \alpha>0, \beta>0, \gamma> 0, \alpha+\beta+\gamma={\pi} \}.
$$

 Let us calculate the invariant area form $\bar \omega$ on this shape space, the result of the symplectic reduction.

\begin{proposition} \label{prop:form}
One has
$$
\bar \omega = \frac{d\alpha \wedge d\beta+d\beta\wedge d\gamma + d\gamma\wedge d\alpha}{6\sin\alpha \sin\beta \sin\gamma}.
$$
\end{proposition} 

\begin{figure}[ht]
\centering
\includegraphics[width=.35\textwidth]{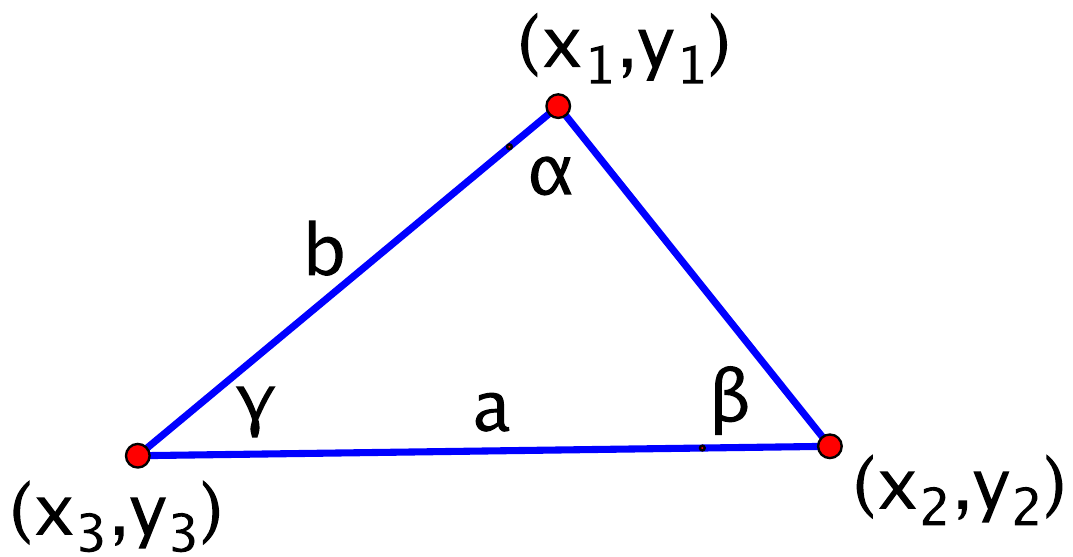}
\caption{To proof of Proposition \ref{prop:form}.}	
\label{tri}
\end{figure}

\proof
Consider a triangle of area 1/2 and, acting by isometries, place its vertex $P_3$ at the origin and the vertex $P_2$ on the positive horizontal axis, see Figure \ref{tri}. Then $ab\sin\gamma=1$ (the area condition), and 
$
x_1=b\cos\gamma,\ y_1=b\sin\gamma,\ x_2=a,\ y_2=x_3=y_3=0.
$
We need to express 
$$
\omega=dx_1\wedge dx_2+dx_2\wedge dx_3+dx_3\wedge dx_1+dy_1\wedge dy_2+dy_2\wedge dy_3+dy_3\wedge dy_1
$$
in terms of the angles $\alpha,\beta,\gamma$.

As an intermediate result, substitute the values of the coordinates of the vertices to the formula for $\omega$ to obtain 
$\omega=da\wedge (\cos\gamma\ db - b \sin\gamma\ d\gamma)$. 

Next, combine the area condition $ab\sin\gamma=1$ with the Sine Rule $a\sin\beta=b\sin\alpha$. This gives
$$
a=t\sin\alpha,\ b=t\sin\beta,\ \ {\rm where}\ \ t^2 \sin\alpha \sin\beta \sin\gamma =1.
$$
Hence
$$
da = \sin\alpha\ dt + t\cos\alpha\ d\alpha,\  db = \sin\beta\ dt + t\cos\beta\ d\beta,\ dt = -\frac{t}{a} \left(  \cot\alpha\ d\alpha + \cot\beta\ d\beta + \cot\gamma\ d\gamma \right)
$$
(the latter is obtained by taking the logarithmic derivative). 
Substitute this to the above intermediate form  and use the fact that $d\gamma = -d\alpha - d\beta$ to obtain, after collecting terms:
$$
\bar \omega=\frac{d\alpha \wedge d\beta}{2\sin\alpha \sin\beta \sin\gamma},
$$
which proves the stated result.
\proofend

In particular, the total area of the shape space is infinite.

Note that the symmetric group $S_3$  acts on the shape space by permutations of the angles. Even permutations commute with the map $\T$, and odd ones conjugate $\T$ with its inverse.

The field $\xi$ descends on the shape space.

\begin{lemma} \label{lm:shape}
The field $\xi$ has an integral 
$$
I:= \frac{{\mathcal P}^2}{4{\mathcal A}} = \cot \frac{\alpha}{2}+\cot\frac{\beta}{2}+\cot\frac{\gamma}{2} = \cot \frac{\alpha}{2}\cot\frac{\beta}{2}\cot\frac{\gamma}{2} .
$$
It  has a unique zero at point $\alpha=\beta=\gamma=\pi/3$, and all its trajectories are closed.
\end{lemma}

\begin{figure}[ht]
\centering
\includegraphics[width=.4\textwidth]{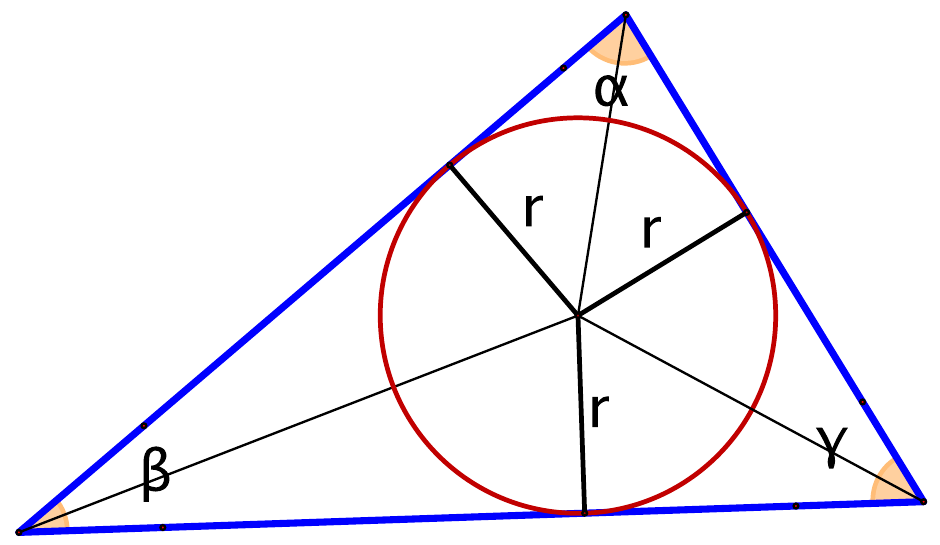}
\caption{Deducing the formula for integral $I$.}	
\label{elem}
\end{figure}

\proof
Let $r$ be the inradius of a triangle. Then 
$$
{\mathcal P} = 2r\left(\cot \frac{\alpha}{2}+\cot\frac{\beta}{2}+\cot\frac{\gamma}{2}\right)
$$
 and ${\mathcal A} = r {\mathcal P} /2$, see Figure \ref{elem}. In particular, the inradius is also conserved by the field $\xi$.
 
 The equality of the sum of cotangents to their product is an easy exercise in trigonometry, and we leave it to the readers.
 
 The function $I$ has a unique critical point $\alpha=\beta=\gamma$ in ${\mathcal S}$, and this is a minimum. All the other level curves of $I$ are closed, and they are the trajectories of the field $\xi$.
\proofend

This lemma implies that, under the flow of $\xi$, every triangle returns to its isometric copy. Furthermore, since the function $I$ is invariant under the $S_3$-action, so are its level curves. This implies that, evolving under the flow of $\xi$, the angles of a triangle undergo all permutations before they return to the original values. 

\begin{figure}[ht]
\centering
\includegraphics[width=.4\textwidth]{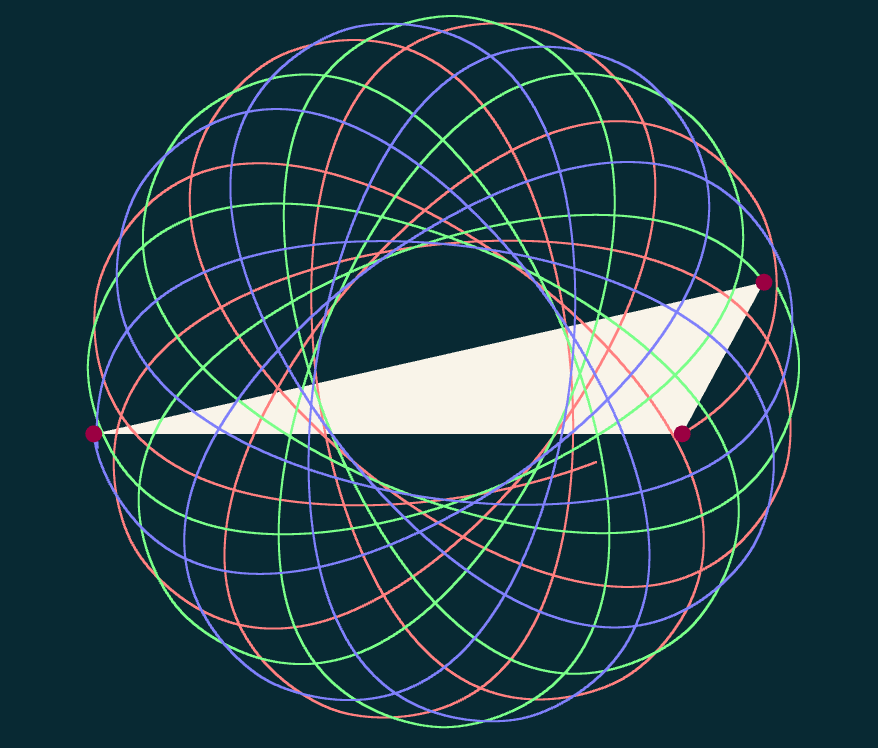}
\caption{The evolution of a triangle under the flow of $\xi$.}	
\label{same}
\end{figure}

This observation explains why the trajectories of the three vertices of a triangle, evolving by the flow of $\xi$, have the dihedral $D_3$-symmetry, see Figure \ref{same}.   

Let us calculate the side lengths of the triangles ${\bf Q}=q_t({\bf P})$ and ${\bf R}=r_t({\bf P})$. 

Denote by $b_i = |Q_{i+1} Q_{i+2}|,\ c_i = |R_{i+1} R_{i+2}|,\ i=1,2,3,$  the side lengths of the triangles ${\bf Q}$ and ${\bf R}$; let $\alpha_i$ be the angles of the triangle ${\bf P}$ and $a_i$ be its side lengths. 

\begin{lemma} \label{lm:QR}
One has
$$
\begin{aligned}
b_i^2=(a_i+t)^2+t^2+t(a_i+t)\frac{a_i^2+a_{i+2}^2-a_{i+1}^2}{a_ia_{i+2}}
=a_i^2+\frac{t(a_i+t)[(a_i+a_{i+2})^2-a_{i+1}^2]}{a_ia_{i+2}},\\
c_i^2=(a_i+t)^2+t^2+t(a_i+t)\frac{a_i^2+a_{i+1}^2-a_{i+2}^2}{a_ia_{i+1}}
=a_i^2+\frac{t(a_i+t)[(a_i+a_{i+1})^2-a_{i+2}^2]}{a_ia_{i+1}}.
\end{aligned}
$$
\end{lemma}

\begin{figure}[ht]
\centering
\includegraphics[width=.4\textwidth]{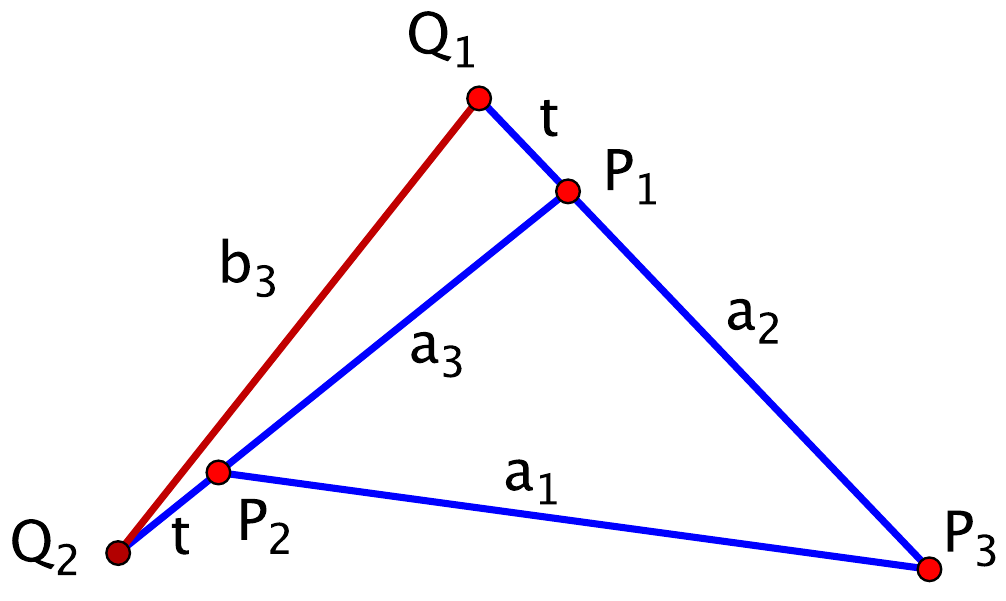}
\caption{To Lemma \ref{lm:QR}.}	
\label{elem2}
\end{figure}

\proof
Let us prove the first formula for $i=3$; the rest is proved in the same way. See Figure \ref{elem2}. 

The Cosine Rule for triangle $Q_1P_1Q_2$ gives 
$$
b_3^2=(a_3+t)^2+t^2+2\cos \alpha_1 t(a_3+t),
$$
and the Cosine Rule for triangle for triangle $P$ gives
$$
a_3^2+a_2^2-2a_3a_2\cos\alpha_1=a_1^2.
$$
Combining the two yields the result.
\proofend

We add to it a formula for the ratio of the area of the triangle ${\bf Q})$ (and ${\bf R}$) to that of ${\bf P}$. We change the notations for the sides from $a_i,\ i=1,2,3$ to  $a,b,c$; accordingly, the angles are denotes by $\alpha,\beta,\gamma$. 

\begin{lemma} \label{lm:rat}
One has
$$
\frac{{\mathcal A}({\bf Q})}{{\mathcal A}({\bf P})} = \frac{abc+t(ab+bc+ca)+t^2(a+b+c)}{abc}=\frac{(t+a)(t+b)(t+c)-t^3}{abc}.
$$
\end{lemma}

\proof
According to Figure \ref{elem2}, one has
\begin{equation} \label{eq:arrat}
{\mathcal A}({\bf Q})={\mathcal A}({\bf P}) + \frac{1}{2} [t(c+t)\sin\alpha+t(a+t)\sin\beta+t(c+t)\sin\gamma].
\end{equation}
One also has
\begin{equation} \label{eq:ar}
{\mathcal A}({\bf P}) = \frac{1}{2} bc\sin\alpha=\frac{1}{2} ca\sin\beta=\frac{1}{2} ab\sin\gamma.
\end{equation}
Thus
$$
\frac{1}{2} \sin\alpha=\frac{{\mathcal A}({\bf P})}{bc},\ \frac{1}{2} \sin\beta=\frac{{\mathcal A}({\bf P})}{ca},\ \frac{1}{2} \sin\gamma=\frac{{\mathcal A}({\bf P})}{ab}.
$$
Substitute to (\ref{eq:arrat}) and simplify to obtain the result.
\proofend

In particular, Lemma \ref{lm:rat} again implies that ${\mathcal A}({\bf Q})={\mathcal A}({\bf R})$.

Consider the $t\to\infty$ limit of the relation $\T$. To do this, we divide the formulas of Lemma \ref{lm:QR} by $t^2$, set $\eps=1/t$, and calculate modulo $\eps^2$. One obtains
$$
b_i^2= \frac{(a_i+a_{i+2})^2-a_{i+1}^2}{a_ia_{i+2}} (1+\eps a_i),\ c_i^2= \frac{(a_i+a_{i+1})^2-a_{i+2}^2}{a_ia_{i+1}} (1+\eps a_i).
$$
This has a curious consequence: {\it in this approximation, the product of the side lengths is a conserved quantity: $\prod b_i = \prod c_i$.} 

The area of a triangle with the sides $a,b,c$ and angles $\alpha,\beta,\gamma$  is given by (\ref{eq:ar}).
Since  area is also a conserved quantity, 
$$
\frac{8{\mathcal A}^3}{(abc)^2} = \sin\alpha \sin\beta \sin\gamma
$$
is, in the $t\to\infty$ approximation, an integral of the resulting vector field on the shape space.







\paragraph{Conics.}
Let ${\bf Q} \T {\bf R}$. One has the following configuration result.

\begin{figure}[ht]
\centering
\includegraphics[width=.63\textwidth]{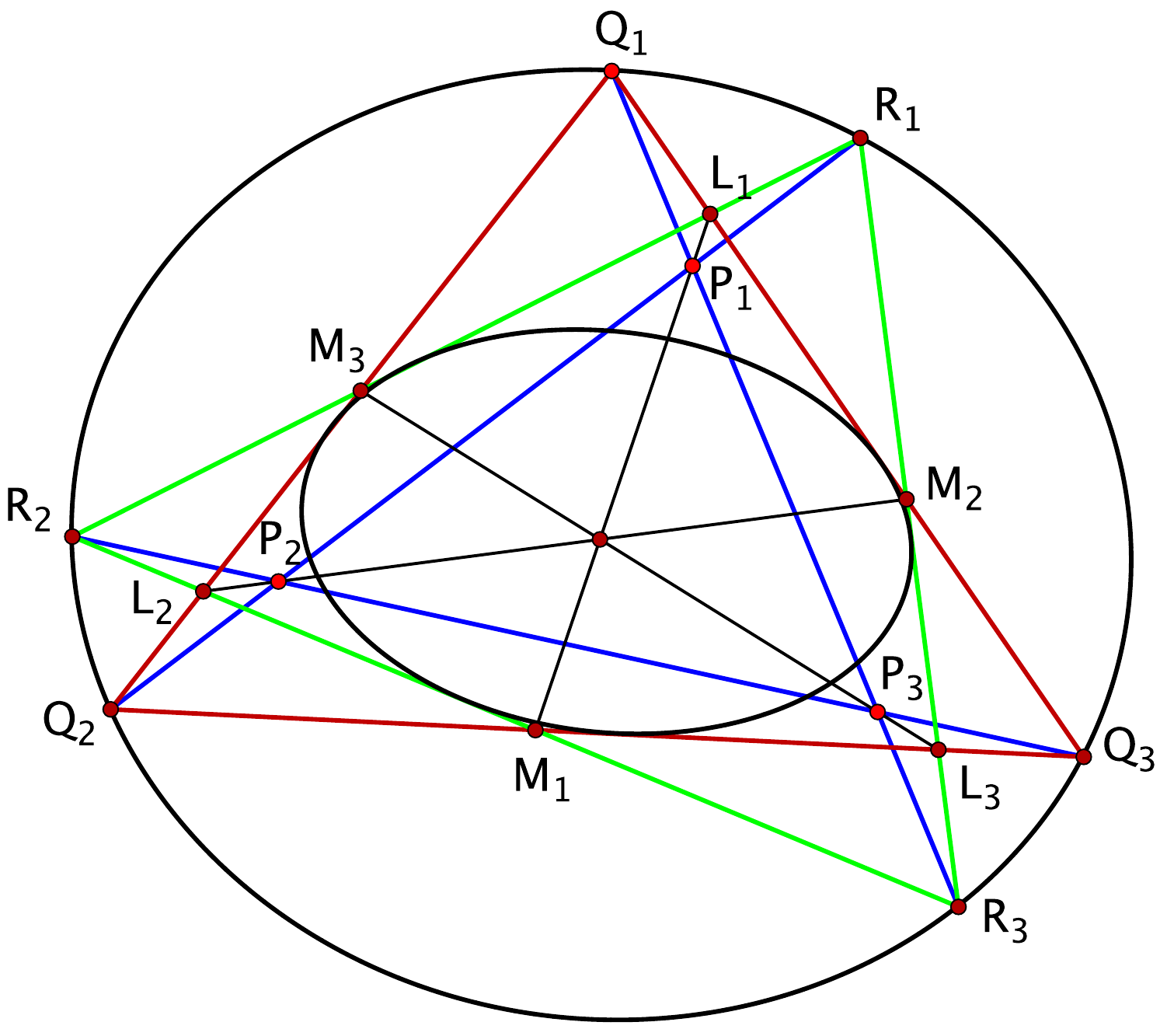}
\caption{To Theorem \ref{conics}.}	
\label{conics}
\end{figure}

\begin{theorem} \label{thm:conics}
The triangles ${\bf Q}$ and ${\bf R}$ are inscribed into a conic and circumscribed about a conic. 
The triples of points $L_1 P_1 M_3, L_2 P_2 M_1$, and $L_3 P_3 M_2$ are collinear, and the
lines $L_1 P_1, L_2 P_2$, and $L_3 P_3$ are concurrent, see Figure \ref{conics}.
\end{theorem}

\proof Recall the Carnot theorem: {\it if a pair of points is chosen on each side of a triangle, then these six points lie on a conic if and only if the product of the three cross-ratios of the quadruples of points on each side equals 1.} Here one uses the following formula of cross-ratio $[a,b,c,d]=(a-c)(c-d)/(a-b)(b-d)$ (out of the six possible choices). See, e.g., \cite{RG} for this and other projective theorems concerning conics.

In our situation, these quadruples are $R_1 P_1 P_2 Q_2, R_2 P_2 P_3 Q_3$, and $R_3 P_3 P_1 Q_1$. Hence
each cross-ratio equals 1, and we conclude that the vertices of the triangles ${\bf Q}$ and ${\bf R}$ lie on a  conic. 

Consider the conic tangent to the five lines $Q_1Q_2, Q_2 Q_3, Q_3 Q_1, R_1 R_2$, and $R_2 R_3$. The triangle ${\bf Q}$ is Poncelet: it is inscribed into a conic and circumscribed about a conic. By the Poncelet porism, point $R_1$ is also a vertex of a Poncelet triangle, hence the line $R_3 R_1$ is tangent to the same inner conic. 

The hexagon $R_1R_2R_3Q_1Q_3Q_2$ is inscribed in a conic. By the Pascal theorem, the points $L_1 P_1 M_1$ are collinear. Likewise with the other two triples of points.

The hexagon $L_1 M_3 L_2 M_1 L_3 M_2$ is circumscribed about a conic. The Brianchon theorem implies that the lines $L_1 M_1, L_2 M_2$ and $L_3 M_3$ are concurrent. 
\proofend

In fact, there are more conics and incidences involved is this picture. We do not dwell on it; see \cite{Ba}.

\section{Even-gons} \label{sect:even}

In this section, $n$ is even. 
Let ${\bf P}=(P_1, P_2,\ldots,P_n)$ be an $n$-gon. Define a vector 
$$
{\mathcal V}({\bf P}) = \sum_{i=1}^n (-1)^{i-1} P_i.
$$ 
This vector is well defined, that is, does not depend on the choice of the origin, and a cyclic permutation of the vertices preserves it, up to the sign. 
The next lemma shows that this vector is an integral of the relation $\T$.

\begin{lemma} \label{lm:vect}
If ${\bf Q} \T {\bf R}$ then ${\mathcal V} ({\bf Q}) = {\mathcal V} ({\bf R})$. 
\end{lemma}

\proof
In the notation of Lemma \ref{lm:Int}, $Q_i = P_i + t e_i, R_i = P_i - t e_{i+1}$. Hence
$$
{\mathcal V} ({\bf Q}) = {\mathcal V} ({\bf P}) + t \sum (-1)^{i-1} e_i = {\mathcal V} ({\bf P}) + t \sum (-1)^i e_{i+1} =
{\mathcal V} ({\bf P}) - t \sum (-1)^{i-1} e_{i+1} = {\mathcal V} ({\bf R}),
$$
as claimed.
\proofend

This lemma provides two integrals of the relations $\T$. In addition, one still has the perimeter integral of the vector field $\xi$.

\begin{lemma} \label{lm:perev}
The perimeter function ${\mathcal P}$ is invariant under the flow of $\xi$.
\end{lemma}

\proof
This follows from the next general fact. 

\begin{figure}[ht]
\centering
\includegraphics[width=.37\textwidth]{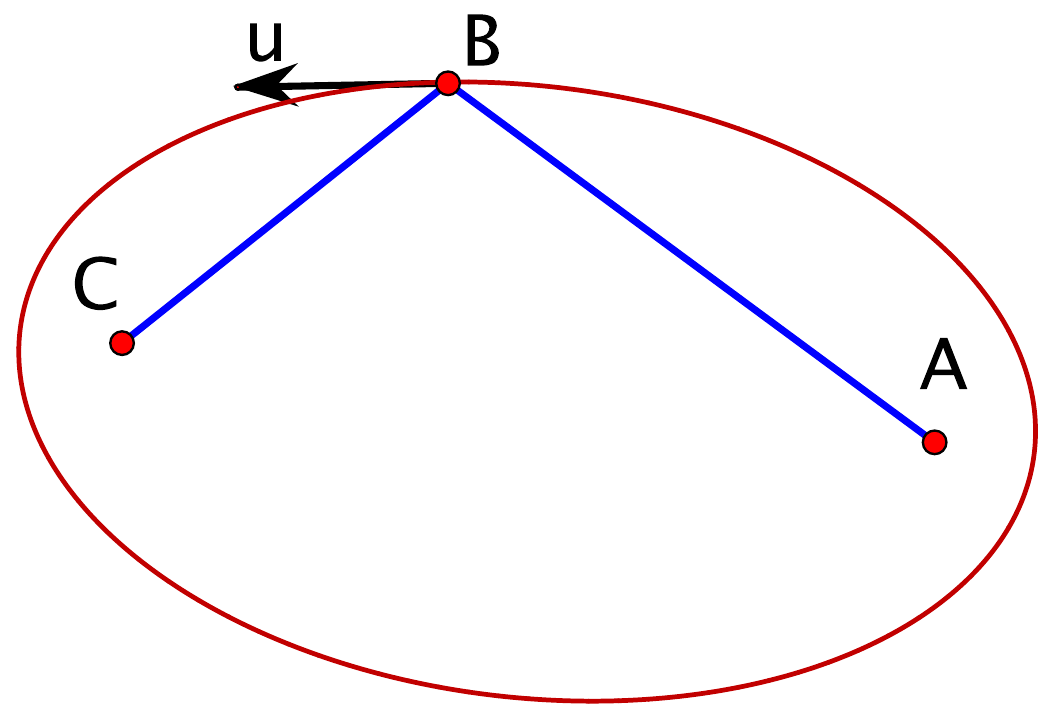}
\caption{To Lemma \ref{lm:perev}.}	
\label{bis}
\end{figure}

Consider a polygonal line $ABC$, let $u$ be a vector in the direction of the bisector of the exterior angle at point $B$, and let $L(B)=|AB|+|BC|$, see Figure \ref{bis}. The level curve of the function $L$ is the ellipse with the foci $A$ and $C$. By the optical property of ellipses, $u$ is tangent to this ellipse, hence $D_u(L)=0$. 

It follows that, since the components of the field $\xi$ bisect the exterior angles of the polygon, its perimeter is invariant under the flow of $\xi$. (This argument is well known in the study of mathematical billiards, see, e.g., \cite{BZ}). 
\proofend

\begin{remark}   
{\rm A necessary and sufficient condition for a polygon ${\bf S}$ to exist such that the vertices of ${\bf P}$ are the midpoints of the sides of ${\bf S}$ is that ${\mathcal V} ({\bf P})=0$. If this condition holds, ${\bf S}$ is defined once one of its vertices is arbitrarily chosen, so there is a 2-parameter family of such polygons ${\bf S}$. Their areas do not depend on the choice involved.
}
\end{remark}


\section{Experimental results, open questions, and conjectures} \label{sect:exp}

In this section we collect open problems and conjectures around complete integrability, along with the results of our experiments. 
The problem of complete integrability has a number of interpretations: it could concern the vector field $\xi$ or the maps $f_t$ in the space of planar $n$-gons or in the moduli spaces ${\mathcal T}$ and ${\mathcal S}$. 

For $n=3$, the space ${\mathcal S}$ is realized as the space of triangles with unit area. This factors out similarities and makes it possible to talk about the value of the parameter $t$.

\paragraph{Odd-gons, space ${\mathcal T}_n$.} As we have seen, this is a symplectic space, the maps $f_t$ and the vector field $\xi$ are symplectic and they preserve the algebraic multi-area ${\mathcal A}$ of a polygon. In addition, $\xi$ is Hamiltonian with the Hamilton function the perimeter ${\mathcal P}$.

\begin{conjecture} \label{conj:intodd}
The maps $f_t$ and the vector field $\xi$ are Liouville integrable, that is,  they possess $n-1$ Poisson commuting integrals, including the algebraic multi-area ${\mathcal A}$ and, in the case of $\xi$, the perimeter function ${\mathcal P}$. 
\end{conjecture}

\begin{figure}[ht]
\centering
\includegraphics[height=.4\textwidth]{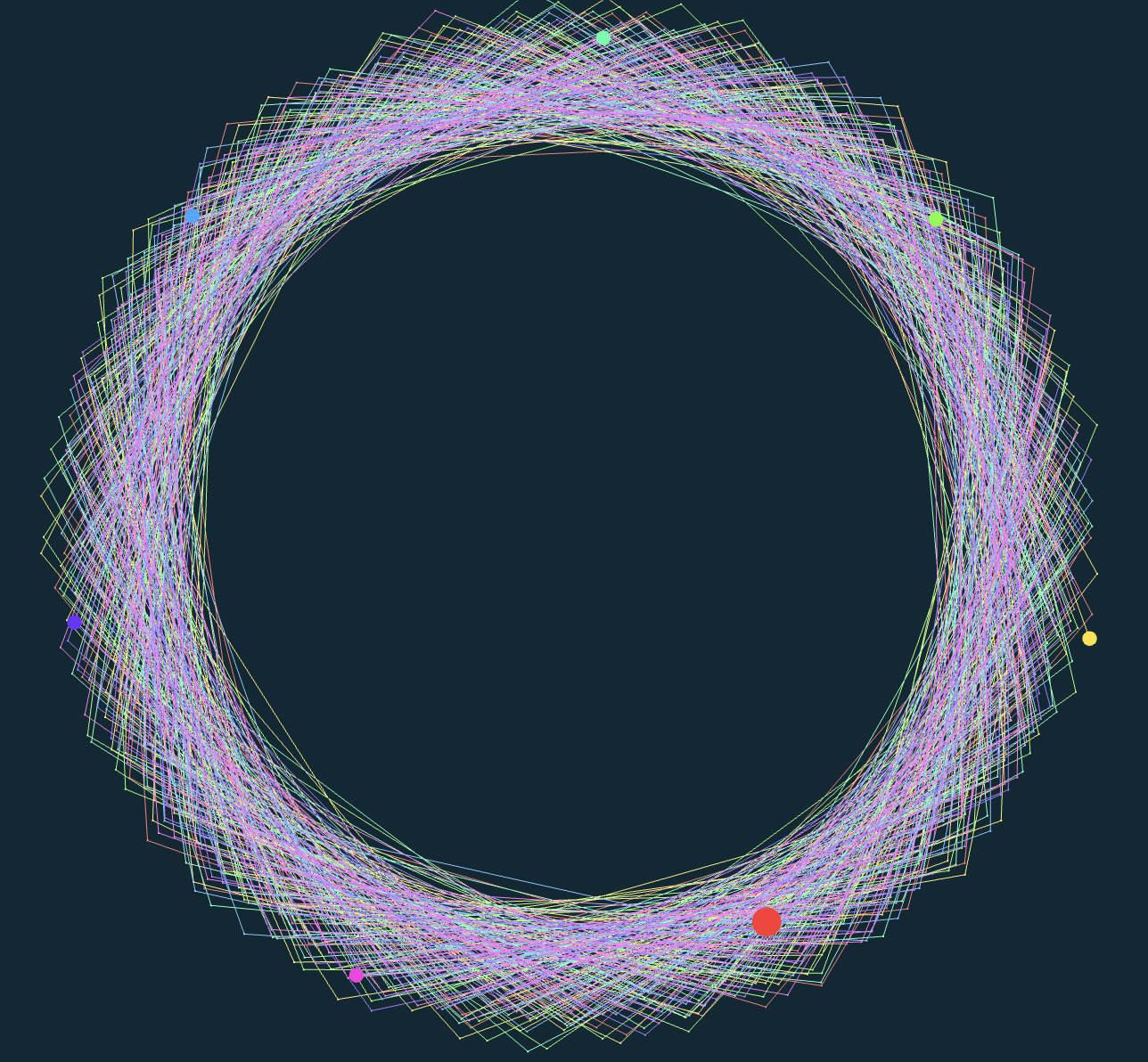}\quad
\includegraphics[height=.4\textwidth]{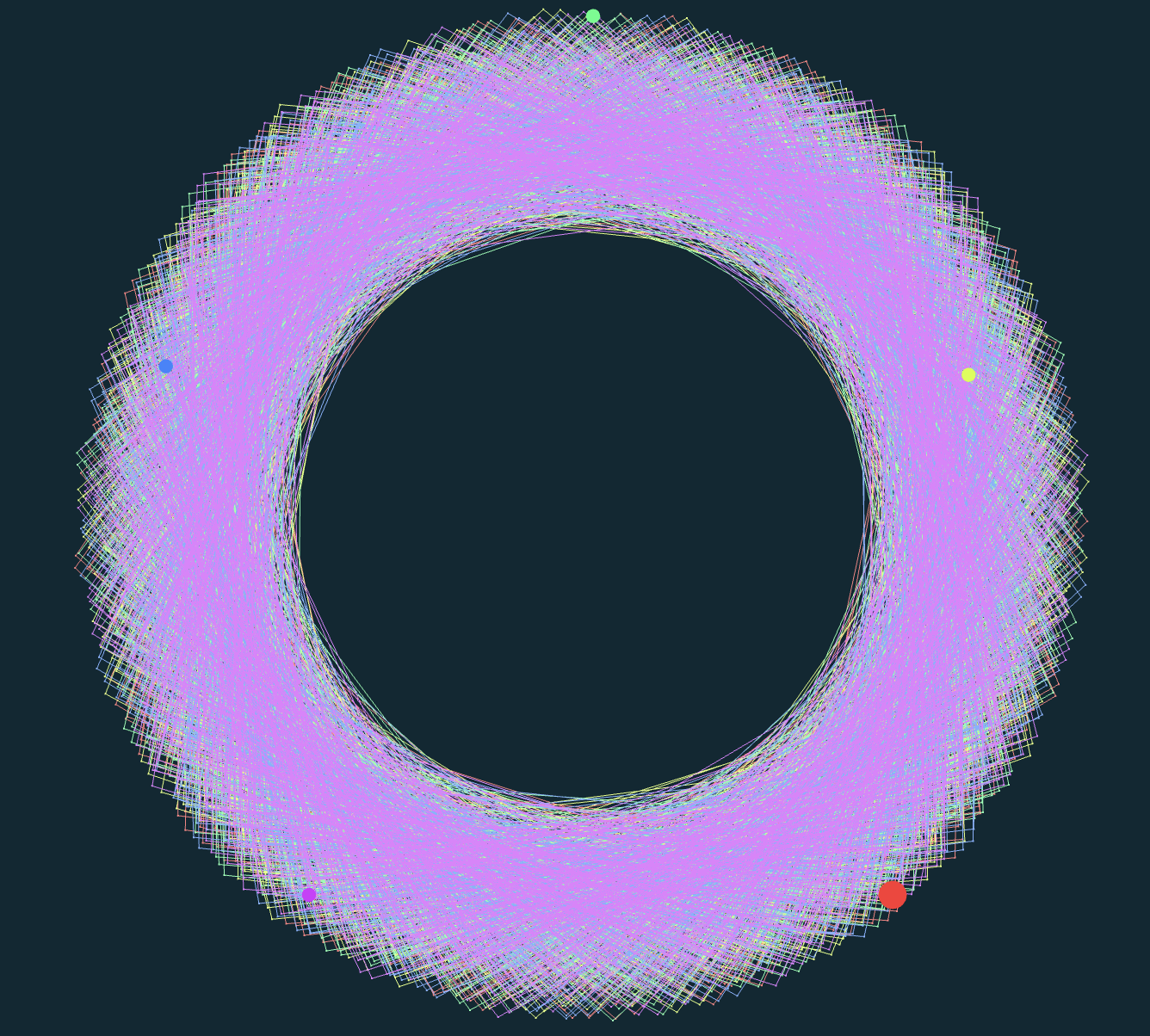}
\caption{$2^{10}$ iterations of a pentagon and of a heptagon.}	
\label{tori}
\end{figure}

If this conjecture holds, the Arnold-Liouville theorem would imply that the motion is quasi-periodic on the common level surfaces of the integrals. This is consistent with our experimental observations, see  Figures \ref{tori}.

\paragraph{Invariant centers.} Figures \ref{same} and  \ref{tori} suggest that the polygons, evolving under the maps $f_t$ and under the flow of the vector field $\xi$, do not drift away, but rather move around certain  invariant centers. 

\begin{conjecture} \label{conj:center}
For every $n$, the maps $f_t$ and the vector field $\xi$ possess two additional integrals, independent of $t$, whose values are the coordinates of the above mentioned centers.
\end{conjecture}


In the case of the vector field $\xi$, we tried to identify these centers in the case of triangles. For this purpose, we turned to the Encyclopedia of Triangle Centers \cite{tri}. We numerically evaluated, with a great precision, the invariant center of a special triangle with the side lengths $6,9,13$. For this triangle, the Encyclopedia of Triangle Centers provides, also with a great precision,  the coordinates of 61,297 centers from this database. We compared our measurement with the data and did not find a single match. 

\paragraph{Failure of  Bianchi permutability.} A version of Bianchi permutability in our setting would be the statement that the maps $f_t$ with different values of the parameter $t$ commute. This is a common feature of 1-parameter families of integrable maps; in particular, it holds for the 
discretization of the bicycle correspondence studied in \cite{TT} and for the families of transformations of polygons studied by us in \cite{AIFT,AFT}. 

This motivated us to investigate Bianchi permutability for the maps $f_t$ in the case of triangles, and we found, numerically, that it did not hold. 

\paragraph{Triangles: when are the maps defined?} 
Given a value of $t$, one wants to describe the set of shapes of triangles for which all the iterations of the map $f_t$ are defined (i.e., such that $t$ remains less that half of every side of the triangles in the orbit). Denote this set by $U_t \subset {\mathcal S}_3$.

For the infinitesimal $t$, that is, for the vector field $\xi$, this not an issue (i.e., $U_t = {\mathcal S}_3$): as we saw in Lemma \ref{lm:shape}, the trajectories of this field are closed curves, and the field is defined for all times.

 \begin{figure}[ht]
\centering
\includegraphics[height=.26\textwidth]{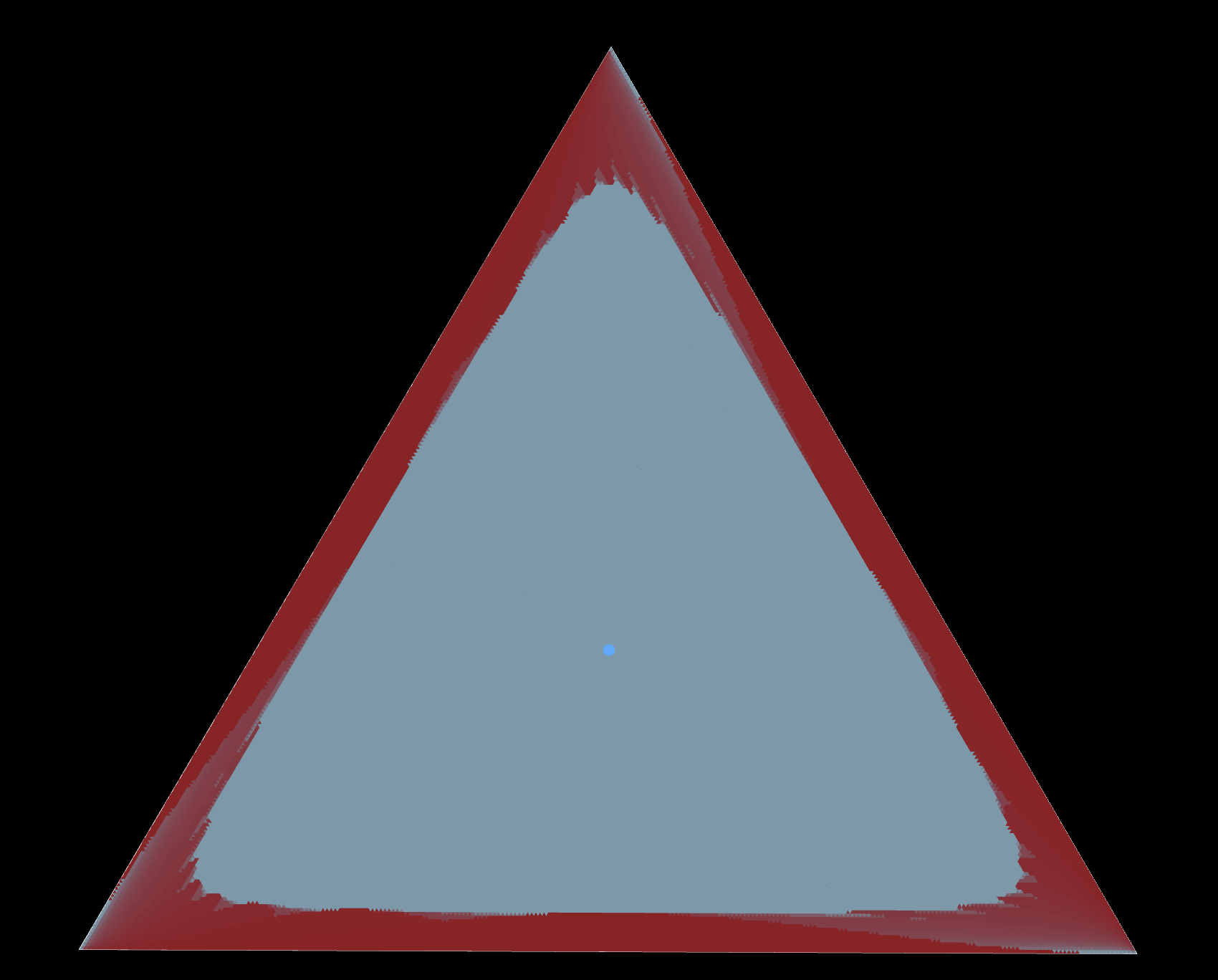}\
\includegraphics[height=.26\textwidth]{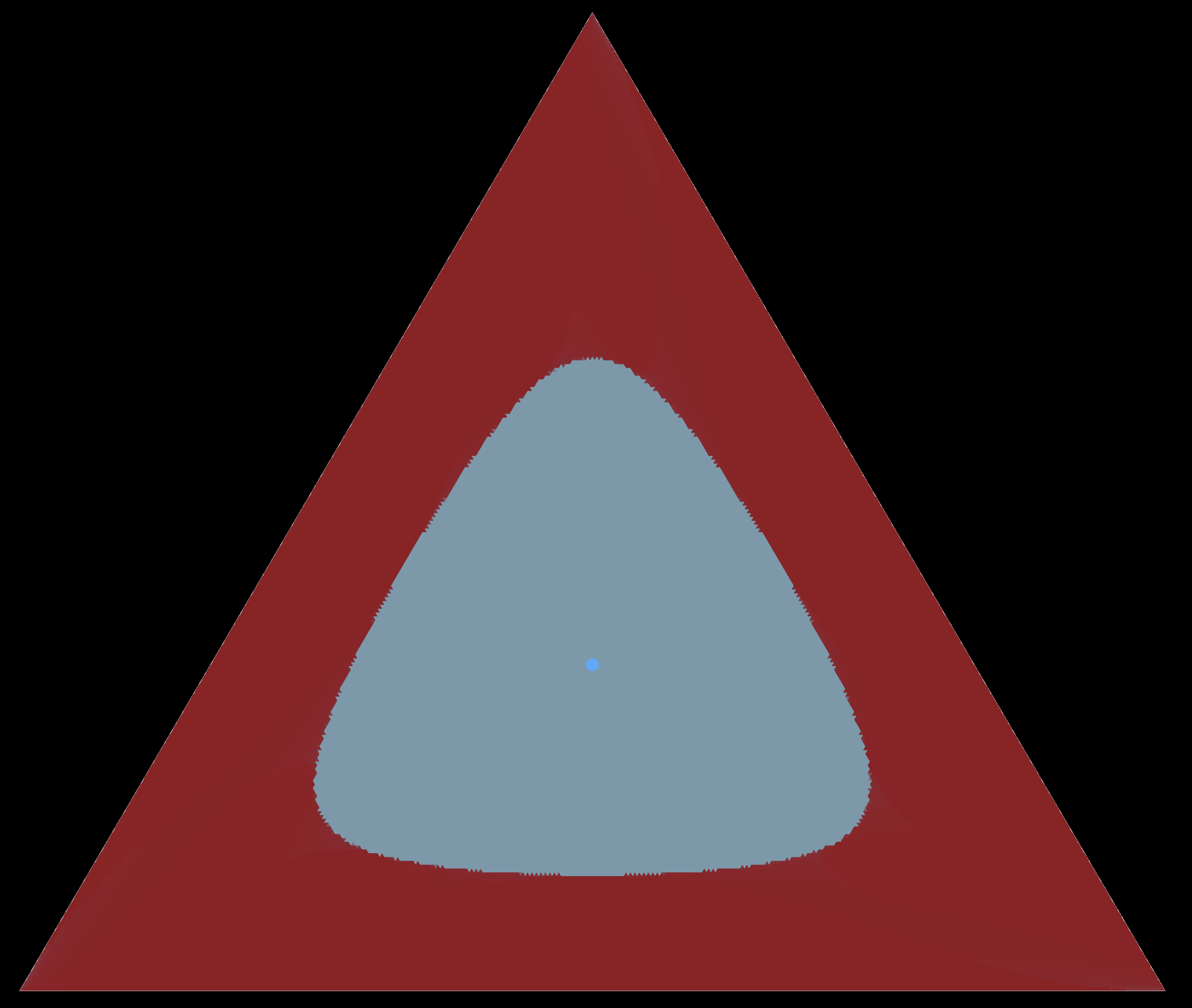}\
\includegraphics[height=.26\textwidth]{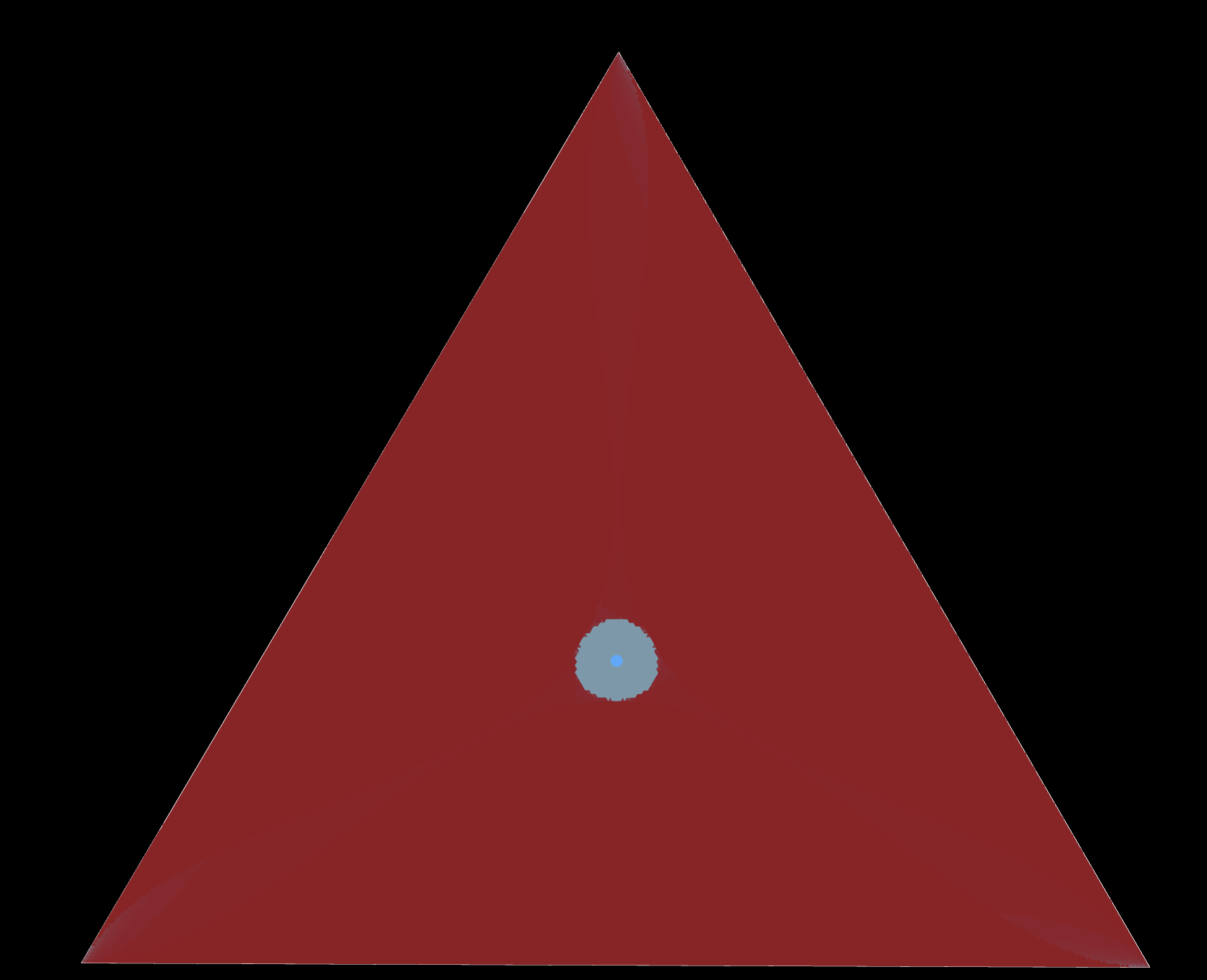}
\caption{The sets $U_t$, in blue, for three increasing values of $t$.}	
\label{zone}
\end{figure}

Figure \ref{zone} shows the sets $U_t$ for various values of $t$. These sets are connected (in fact, convex) and, as $t$ increases, they shrink toward the central point corresponding to the equilateral triangle. 

It would be interesting to describe the curves that bound the sets $U_t$. If Conjecture \ref{conj:intodd} holds, perhaps they are level curves of the integrals of the maps $f_t$.

Of course, one can ask the same question about $n$-gons with $n>3$. We have nothing to say about this intriguing problem.

\begin{figure}[ht]
\centering
\includegraphics[height=.4\textwidth]{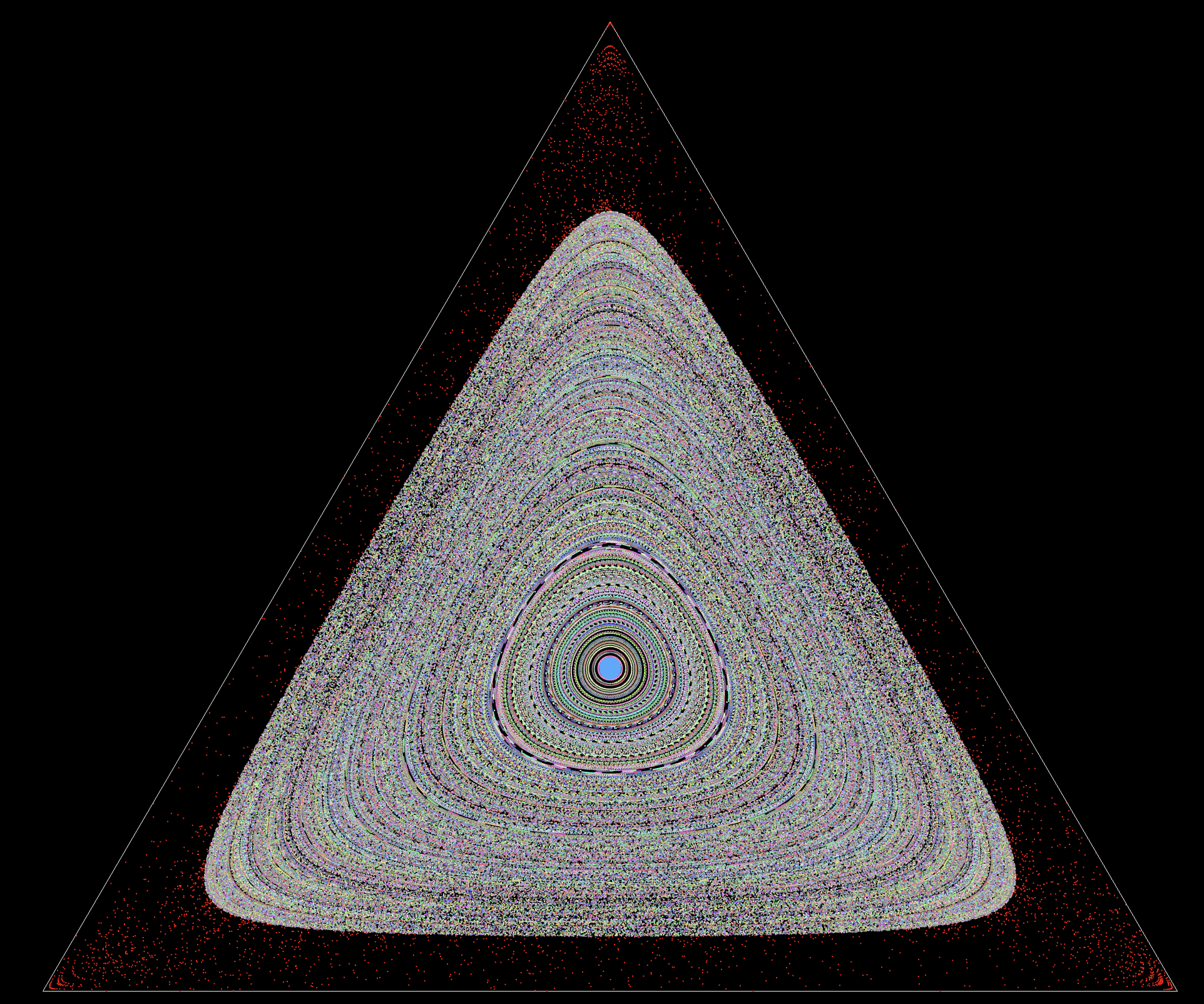}\quad
\includegraphics[height=.4\textwidth]{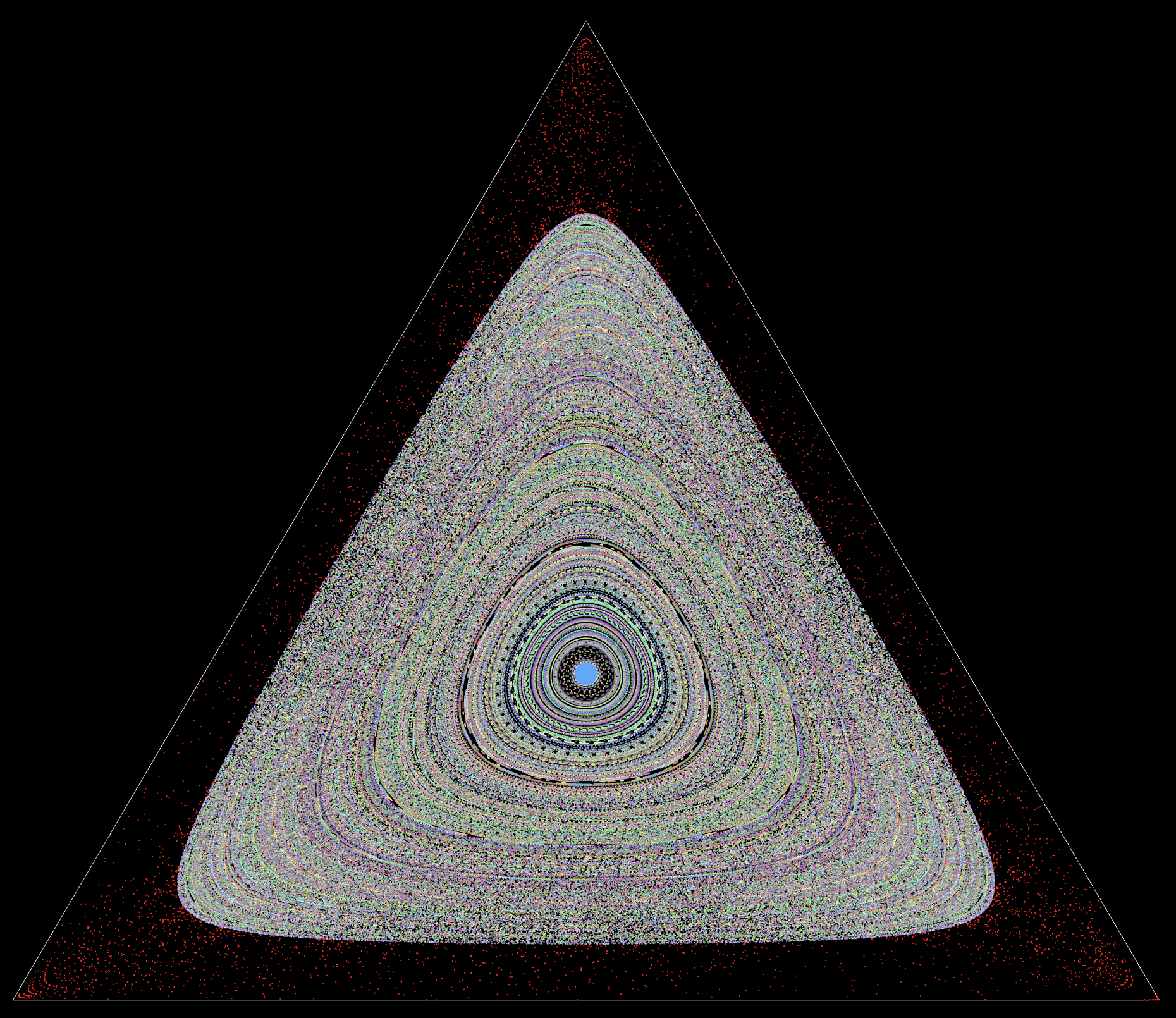}
\caption{Phase portraits of the maps $f_t: {\mathcal S}\to {\mathcal S}$ with two different values of $t$.}	
\label{phase}
\end{figure}

\paragraph{Triangles: phase portraits.} Figure \ref{phase} shows phase portraits of the maps $f_t$ for different values of $t$.
We formulate a particular case of Conjecture \ref{conj:intodd} as a separate conjecture.

\begin{conjecture} \label{conj:tri}
The map $f_t: {\mathcal S}_3 \to {\mathcal S}_3$ has a second integral (depending on $t$), commuting with the area integral ${\mathcal A}$.
\end{conjecture} 


\paragraph{Triangles: periodic orbits.}

\begin{figure}[ht]
\centering
\includegraphics[width=.3\textwidth]{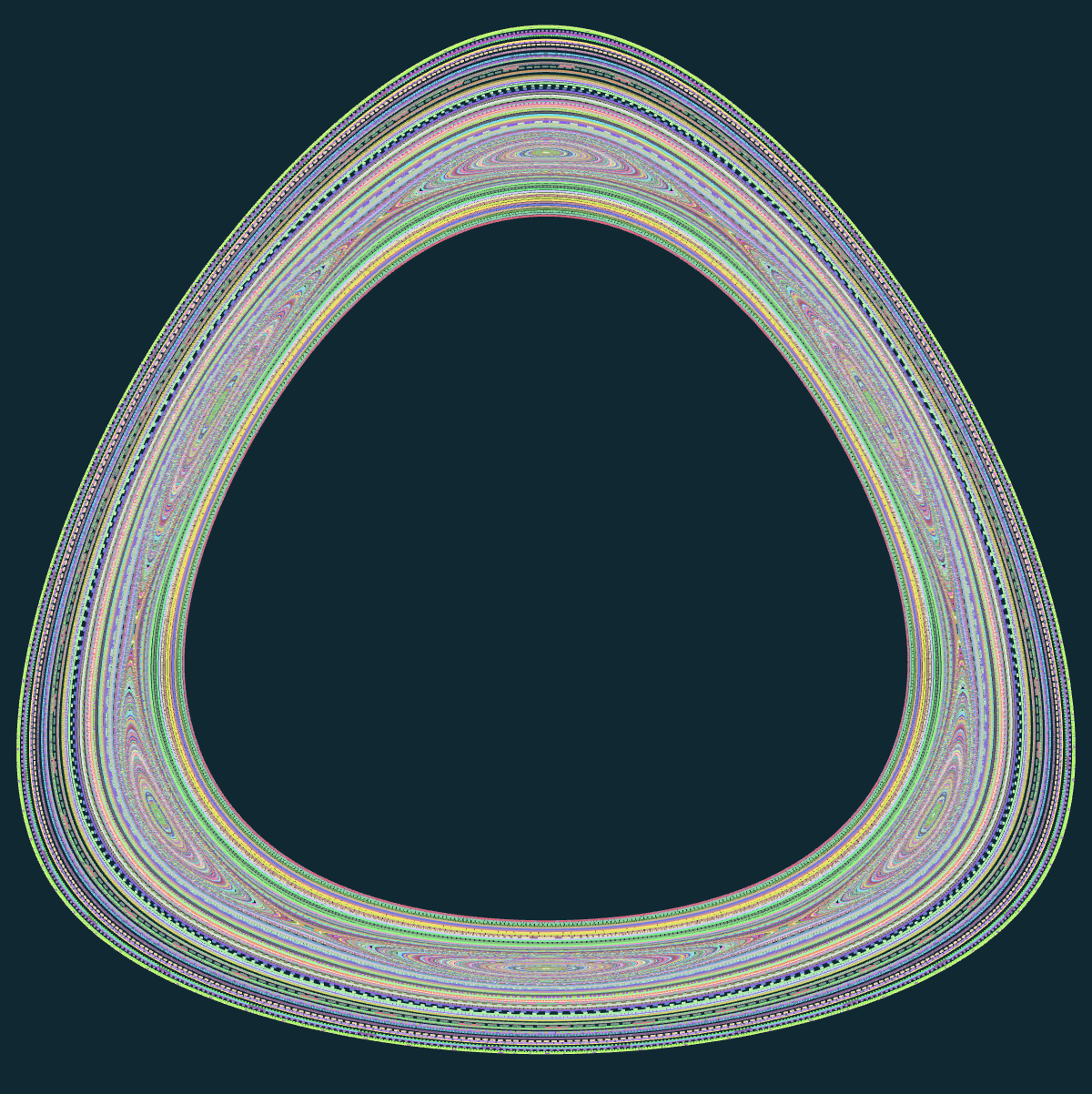}\quad
\includegraphics[width=.3\textwidth]{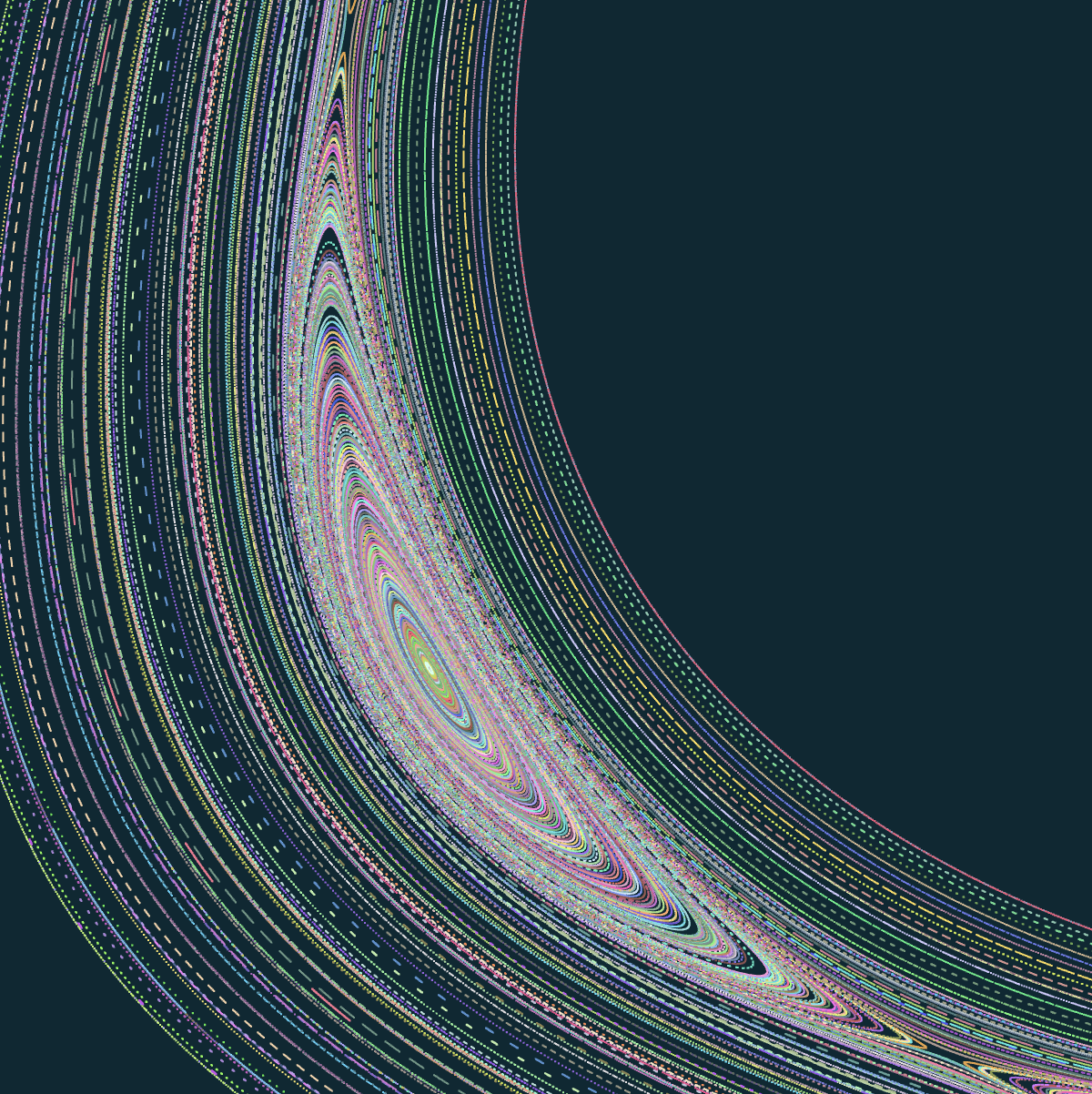}\quad
\includegraphics[width=.3\textwidth]{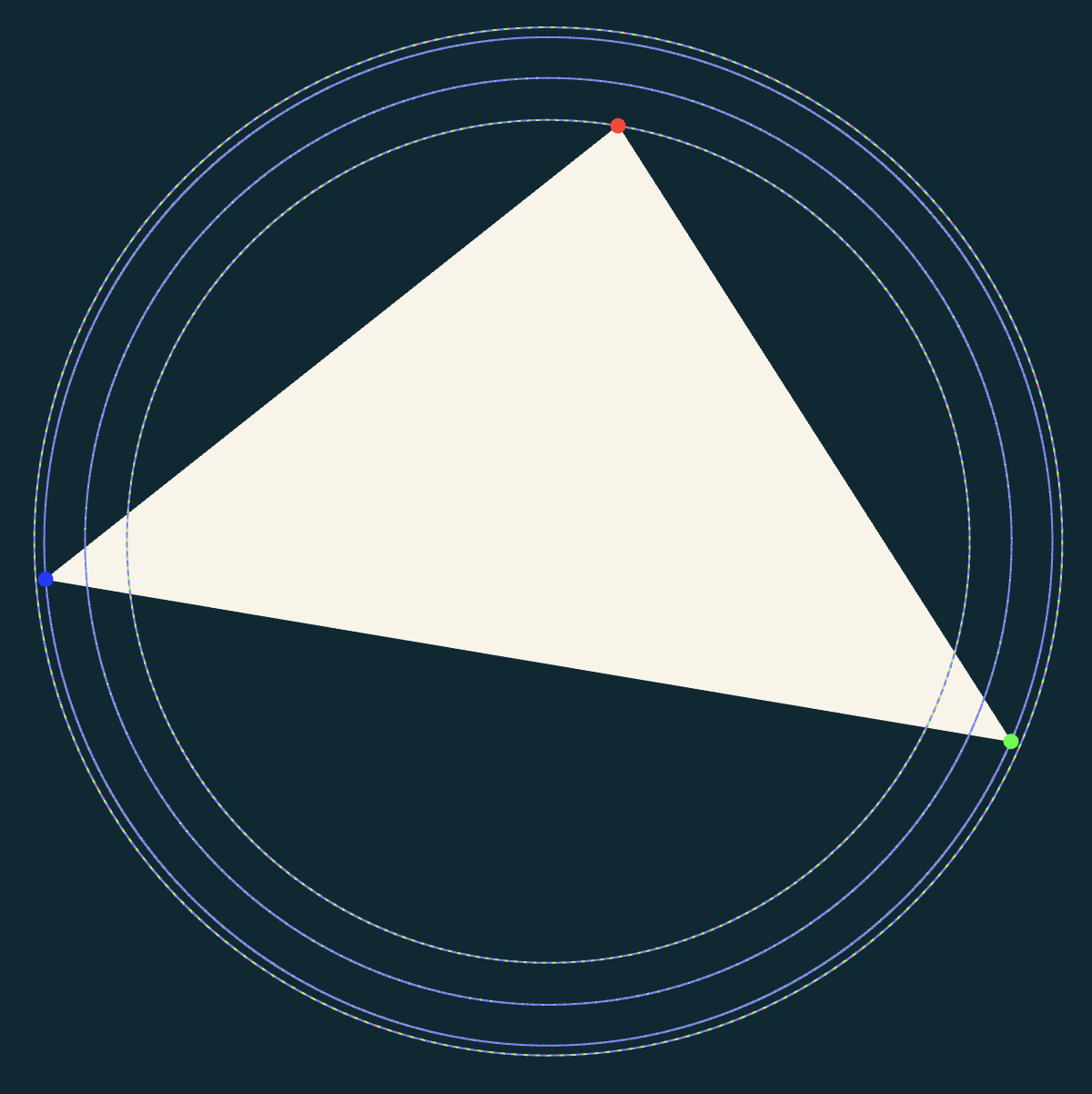}

\caption{
	Partial phase portrait and detail for $t=0.56$ showing islands about an
	apparent 6-periodic orbit. Only orbits reasonably close to the periodic one
	are plotted. The right-most image shows the periodic orbit plotted in the plane.
}	
\label{fig:period-6}
\end{figure}

Motivated by phase portraits like the one shown in Figure \ref{fig:period-6}, we investigated the prospect of a 6-periodic orbit in
$\mathcal S_3$ in which all images are isosceles triangles of unit area. By the symmetry of
$\mathcal{S}_3$, we see that such an orbit consists of exactly two isosceles
triangles up to vertex relabeling; let us call them $\bf Q$ and $\bf R$. Denote
$q_t^{-1}(\bf{Q})=\bf P$ and $r_t^{-1}(\bf{Q})=\bf{P'}$. Since $\bf Q$ is isosceles, we
have that $\bf P$ and $\bf P'$ must be
mirror images of one another along the axis of symmetry of $\bf Q$. 

One way to phrase the problem analytically is that we are looking for
$(a_1,a_2,a_3,t)$ (where $a_i$ are the side lengths of $\bf P$) such that
\begin{enumerate}
	\item $q_t(\bf{P})$ is isosceles:
		$a_1^2+\frac{t(a_1+t)[(a_1+a_3)^3-a_2^2]}{a_1a_3}
		=a_2^2+\frac{t(a_2+t)[(a_2+a_1)^3-a_3^2]}{a_2a_1}$
	\item $r_t(\bf{P})$ is also isosceles:
		$a_1^2+\frac{t(a_1+t)[(a_1+a_2)^3-a_3^2]}{a_1a_2}
		=a_3^2+\frac{t(a_3+t)[(a_3+a_1)^3-a_2^2]}{a_3a_1}$
	\item $q_t(\bf{P})$ has normalized area:
		$\frac{b_3}{2}\sqrt{b_1^2-\frac{b_3^2}{4}}=1$
\end{enumerate}
where we have used the notation and result of Lemma \ref{lm:QR}. Three relations
on four unknowns suggests a single degree of freedom, which matches the phase
portrait experiments we performed. Although we
were unable to solve these equations analytically, we were able to find
numerical solutions with Mathematica. We obtained evidence for the following conjecture:

\begin{conjecture}
	There exist values $t_m\approx0.53$ and $t_M\approx0.57$ such that for every
	$t\in(t_m,t_M)$, there exists a unique non-equilateral triangle ${\bf{P}}(t)$ such that
	$q_t(\bf{P})$ and $r_t(\bf{P})$ are both isosceles and have unit area.
\end{conjecture}

Experimentally, as $t$ approaches $t_m$, ${\bf{P}}(t)$ approaches the equilateral
fixed point. Conversely, if $t>t_M$, the condition that all side lengths of the
triangles in the orbit exceed $2t$ is violated.

We may of course plot these orbits in the plane as well as in the shape space
(see the third image in Figure \ref{fig:period-6}).
In the plane, we notice that the vertices of the iterates lie on four
concentric circles. This numerical result is consistent with the notion that the
squared map acts on each of $\bf Q$ and $\bf R$ by rotation about an invariant
circle. This observation therefore serves as experimental evidence for
Conjecture \ref{conj:center}.



\paragraph{Quadrilaterals: symmetry.}

\begin{figure}[ht]
\centering
\includegraphics[width=.3\textwidth]{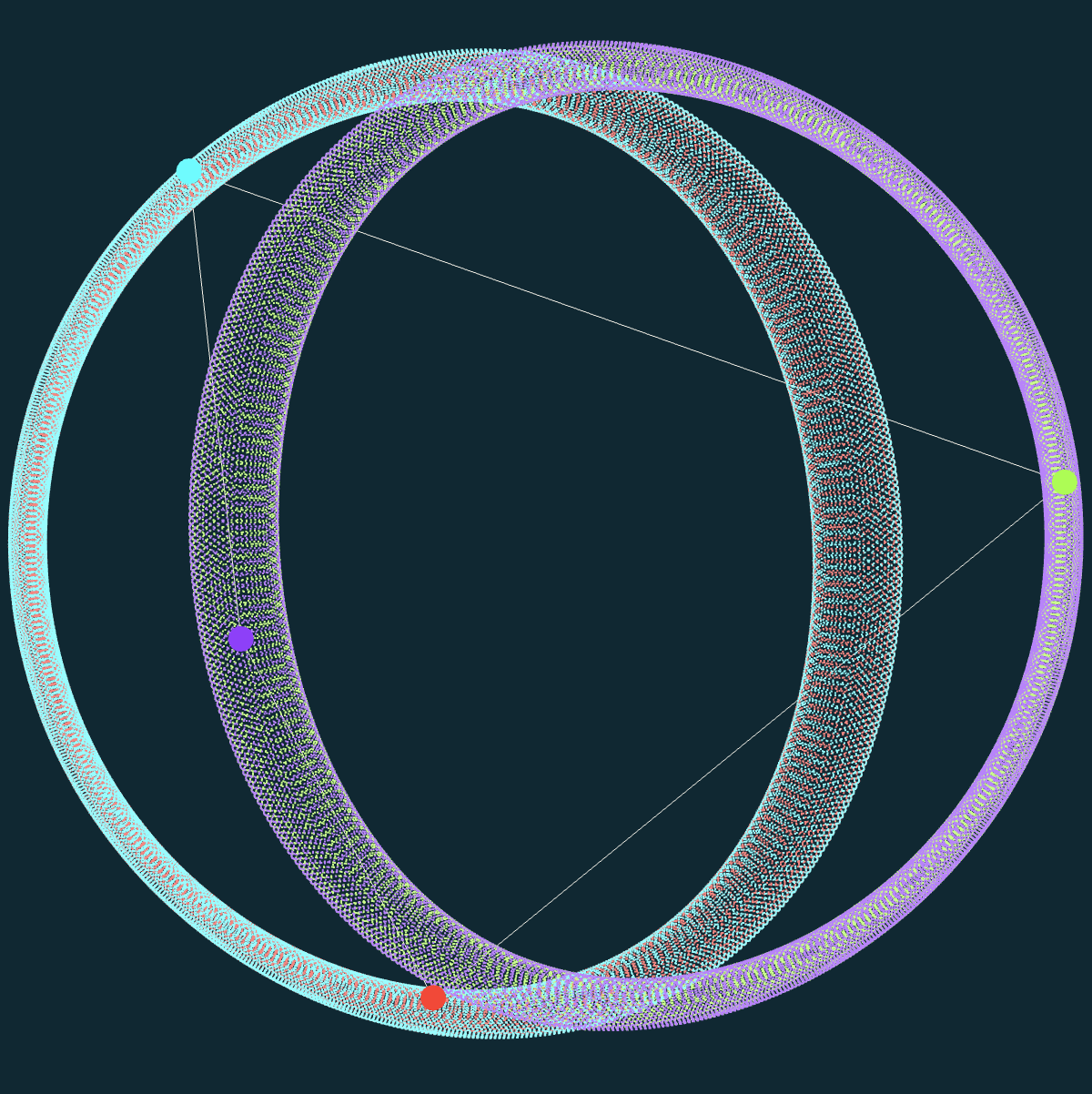}\quad
\includegraphics[width=.3\textwidth]{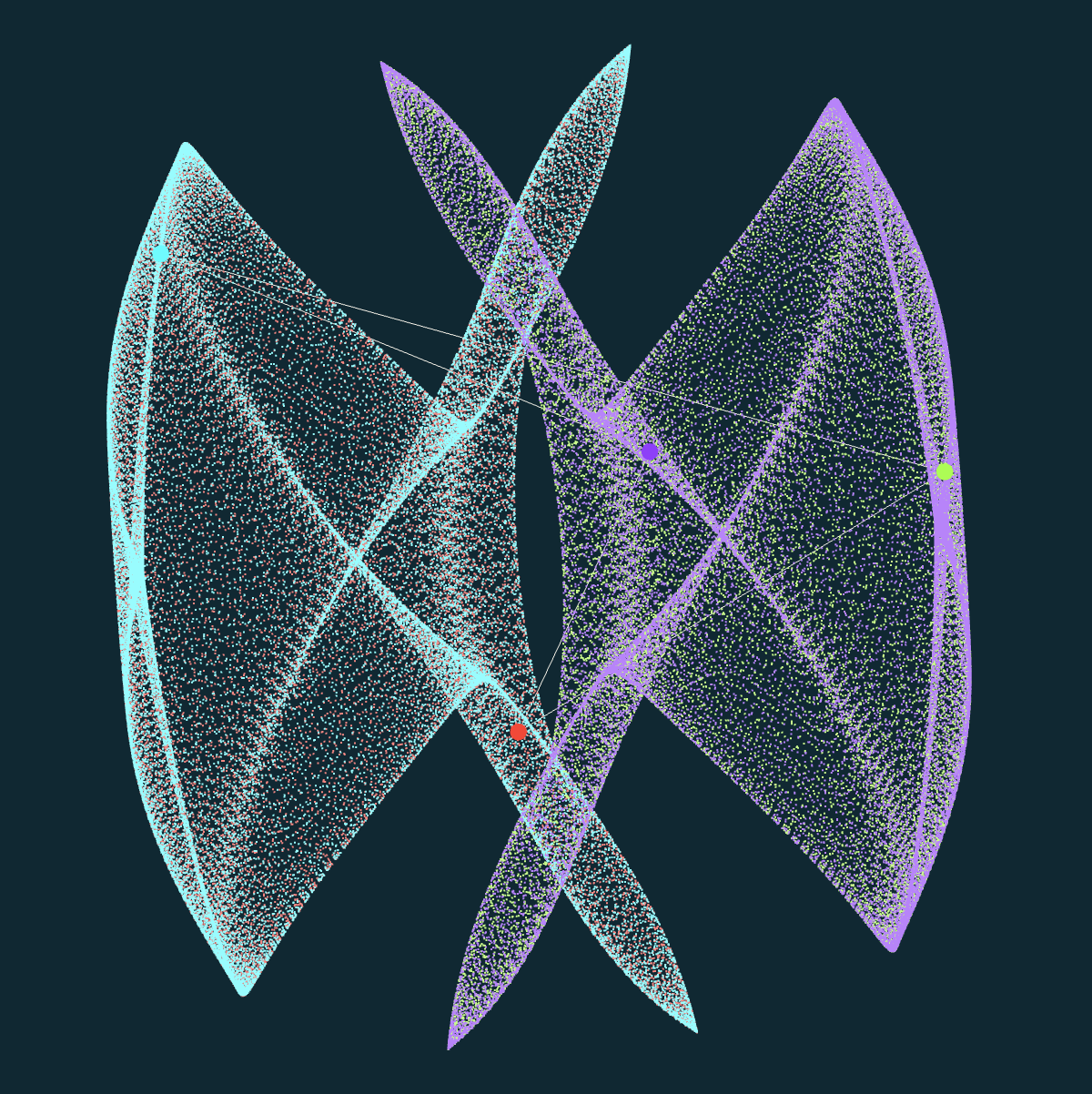}\quad
\includegraphics[width=.3\textwidth]{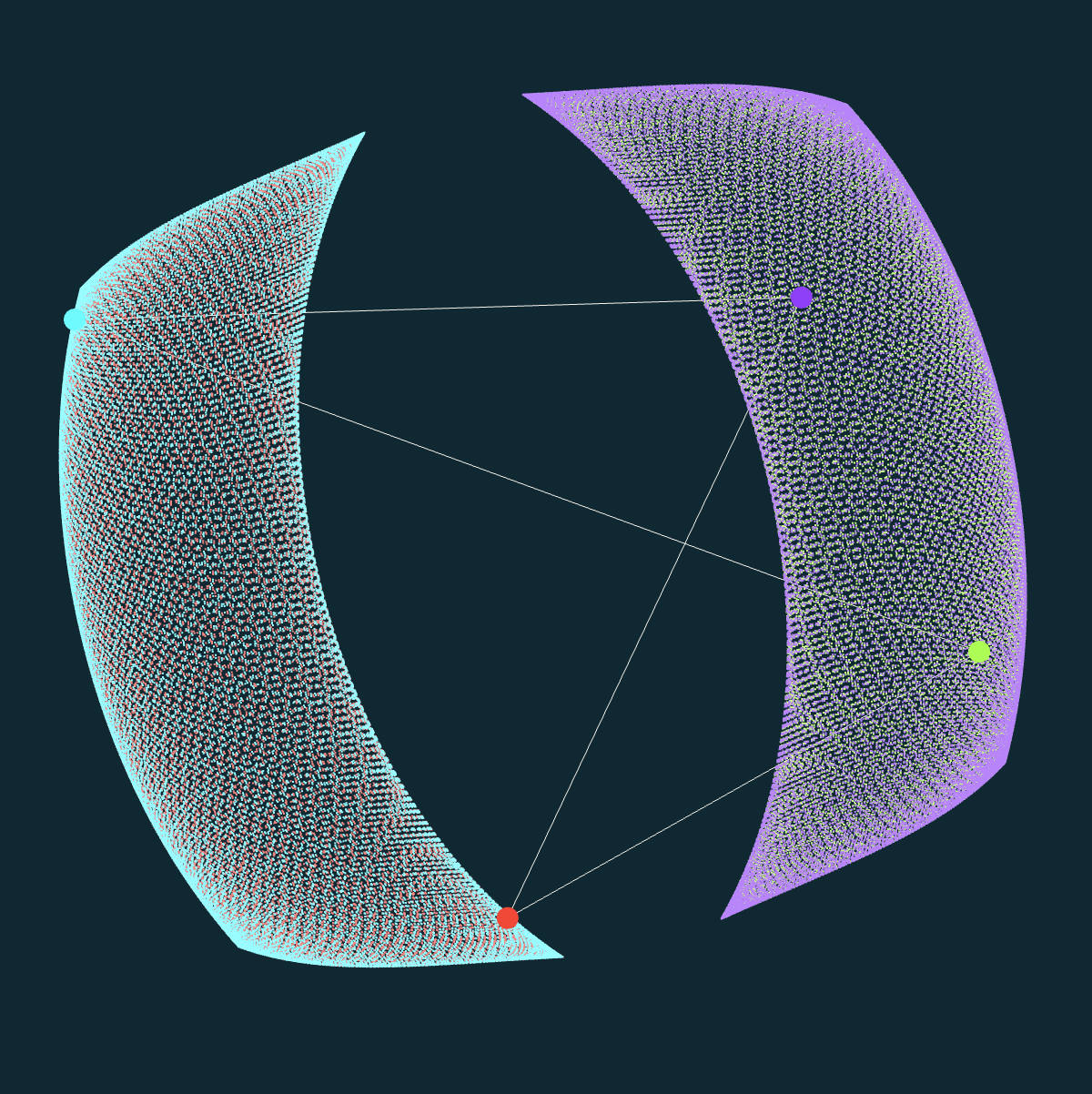}
\caption{
	Orbits of several quadrilaterals under maps $f_t$. The original
	quadrilateral is shown with larger dots for the vertices and white edges,
	and the vertices of the first $2^{15}$ images are plotted.
}	
\label{fig:quads}
\end{figure}

Experiments with the maps $f_t$ applied to quadrilaterals produce orbits
possessing some apparent symmetry. Each orbit (three examples are shown in Figure
\ref{fig:quads}) has
two axes of mirror symmetry which conjecturally intersect at a conserved center
and appears to be made of two identical components, one rotated by
$180^\circ$  relative to the other. One component is made up of images of the
first and third vertices; the other of images of the second and fourth vertices.

In general, we have very little to say, either theoretically or experimentally, about the dynamics of even-gons. This is also an intriguing open problem.

\end{document}